\documentclass{amsart}
\usepackage{graphicx}
\usepackage{amssymb}
\usepackage{amsmath}
\usepackage{amsthm}
\usepackage{enumerate}
\usepackage{bbm}
\usepackage{caption}
\usepackage{subcaption}
\usepackage{calligra}
\usepackage{xcolor}
\usepackage{multicol}
\usepackage{MnSymbol}
\DeclareMathAlphabet{\mathcalligra}{T1}{calligra}{m}{n} 

\newtheorem{thm}{Theorem}[section]

\newtheorem{prop}[thm]{Proposition}
\theoremstyle{definition}
\newtheorem{defn}[thm]{Definition}

\newtheorem*{defn*}{Definition}
\newtheorem*{rems*}{Remarks}
\newtheorem*{rem*}{Remark}

\newtheorem{col}[thm]{Corollary}

\numberwithin{equation}{section}

\newcommand{\Eq}{{\mathrm{E}}}
\newcommand{\E}{\mathbb{E}}

\newcommand{\Css}{{\mathrm{CSS}}}

\usepackage{listings}
\usepackage{xcolor}

\lstdefinelanguage{Mathematica}{
  morekeywords={Plot, Sin, Cos, Exp, Pi, E, I, Table, Do, If, While, For, Print, supportfunction, gamma, tang},
  sensitive=true,
  morecomment=[l]{(*},
  morecomment=[s]{(}{)},
  morestring=[b]",
  morestring=[d]',
  alsoletter={\\, \[ \]},
  alsodigit={., e, E}
}

\begin{document}

\title[Artistic Aspects of the Wigner Caustic and CSS] {Artistic Aspects of the Wigner Caustic and the Centre Symmetry Set}
\author[I. Danielewska, D. Po\l{}awski, D. Sterczewska, and M. Zwierzy\'nski]{Iza Danielewska, Dawid Po\l{}awski, \\ Dominika Sterczewska, and Micha\l{} Zwierzy\'nski$^{\pentagram}$}
\address{Warsaw University of Technology\\
Faculty of Mathematics and Information Science\\
ul. Koszykowa 75\\
00-662 Warsaw, Poland}

\email{Iza.Danielewska.stud@pw.edu.pl, Dawid.Polawski.stud@pw.edu.pl, \linebreak dkaczmarska97@gmail.com, Michal.Zwierzynski@pw.edu.pl}

\thanks{$^{\pentagram}$Corresponding author}

\subjclass[2010]{52A38, 52A40, 58K70}

\keywords{Centre Symmetry Set, envelope, singularities, string art, Wigner caustic}

\begin{abstract}
The Wigner caustic and the Centre Symmetry Set of a closed smooth planar curve are known singular sets which generically admit only cusp singularities. Applications of these objects in semi-classical quantum physics, in chaos theory, in singularity theory, in convex geometry, have been studied since the 1970s until today. These sets can be viewed as envelopes of special families of lines and thanks to that they have many geometric artistic values.
\end{abstract}

\maketitle

\section{Introduction}

The notion of the singular set known as the Wigner caustic was first introduced in the celebrating paper \cite{B1} by Berry on the semi-classical limit of Wigner's phase-space representation of a quantum states. In the semi-classical limit the Wigner function of the classical correspondence $\mathcal{C}$ of a pure quantum state takes on high values at points in a neighborhood of $\mathcal{C}$ and also in a neighborhood of the Wigner caustic of $\mathcal{C}$, which physical properties are studied in many papers since then (see \cite{Chaos} and the literature therein). The Wigner caustic can be viewed as the locus of midpoints of chords connecting parallel pairs on $\mathcal{C}$, i.e. the two points with parallel tangents to $\mathcal{C}$ (\cite{DJRR1, DMR1, DRR1, Romero, DZ-singular, DZ-WC, GH1, GiblinReeve, GiblinReeve2, GZ1}). The Wigner caustic is a special case of an affine $\lambda$-equidistant set (for $\lambda=0.5$) which is the locus of points dividing the same chords in the fixed ratio $\lambda$. The surprising results is that the oriented area of the Wigner caustic of an oval improves the classical isoperimetric inequality known to ancient Greeks for over $2000$ years (\cite{Zwierz1, Zwierz2, Zwierz3}) and improves some other types of isoperimetric inequalities (\cite{CGR, Zhang1, Zhang2}). The properties of the middle hedgehog, which is a generalization of the Wigner caustic, were studied in \cite{Schneider1, Schneider2}. Furthermore, the Wigner caustic gives one of the two constructions of bi-dimensional improper affine spheres, which can be generalized to higher dimensions and are very closely related to solutions of the famous Monge-Amp\'ere equation (\cite{CDR-IAS}). The singular points of all affine equidistants create the Centre Symmetry Set, which is a singular generalization of the notion of the center of symmetry. This set is widely studied in \cite{DR1, JJR1, DZGauss} and in the literature therein.

In this paper, we will define the affine $\lambda$-equidistants, including the Wigner caustic, and the Centre Symmetry Sets as envelopes of some families of affine chords and we will presents some geometric theorems related to these sets and visualization of these theorems. 
Furthermore, we will present artistic pictures based on these sets, which are so important in geometry. Our designs will be based on support functions, which can simplify global the parametrization of described sets.

The tools used to present the artistic aspects of the Wigner caustic and the Centre Symmetry Set are based on the concept of the envelope. The envelope, being a curve that is tangent to each member of a family of curves, is utilized in the artistic technique known as "String Art." String Art involves creating intricate geometric patterns by threading strings between points (often nails) arranged along lines or shapes. By precisely positioning these points and appropriately tensioning the strings, unique and aesthetically appealing patterns are formed, which reflect the mathematical properties of envelopes. This technique not only expresses artistic creativity but also illustrates mathematical phenomena, making it an excellent tool for exploring the artistic aspects of the discussed topics (see \cite{ArtArticle} and the literature therein).

\section{Properties of the Wigner Caustic and the Centre Symmetry Set}

In this section we will describe the definition of affine $\lambda$-equidistants and the Centre Symmetry Set in details, exploring their fundamental characteristics and properties

\begin{figure}[h]
    \centering
    \begin{subfigure}[h]{0.32\textwidth}
        \centering
        \includegraphics[width=\textwidth]{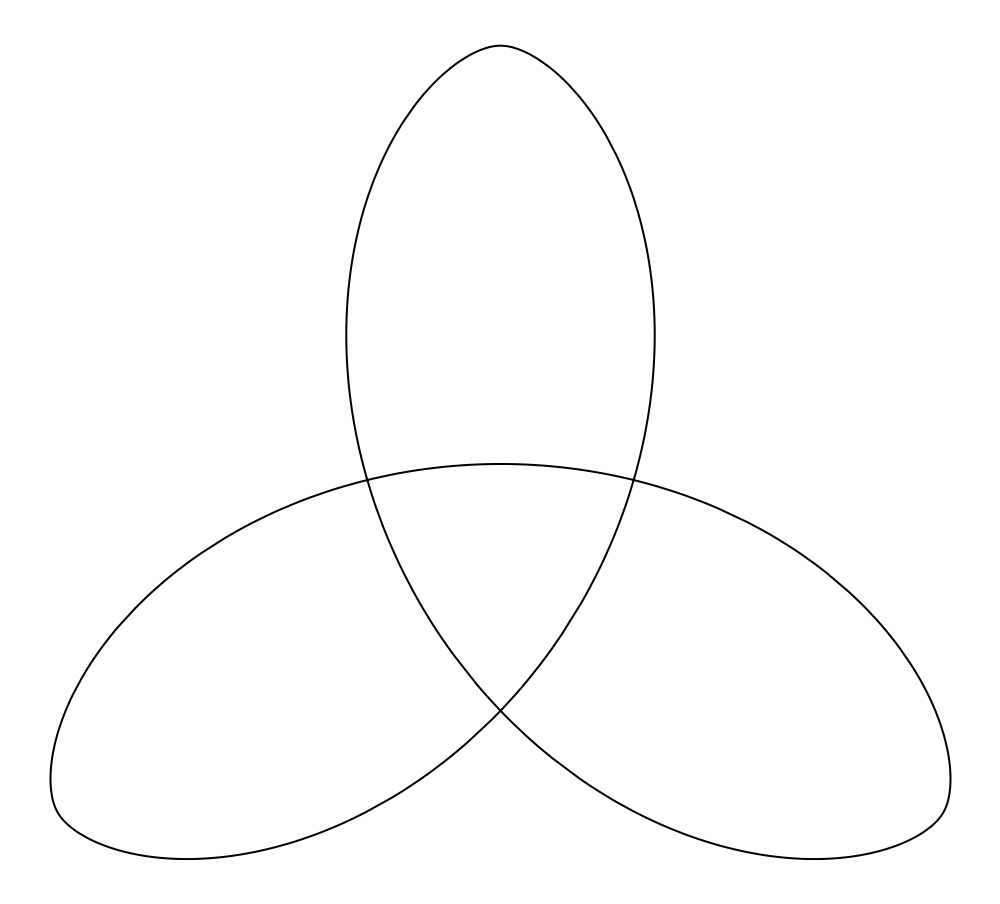}
        \caption{A $2$-rosette}
        \label{fig:Fig_rosettes01}
    \end{subfigure}
    \hfill
    \begin{subfigure}[h]{0.32\textwidth}
        \centering
        \includegraphics[width=\textwidth]{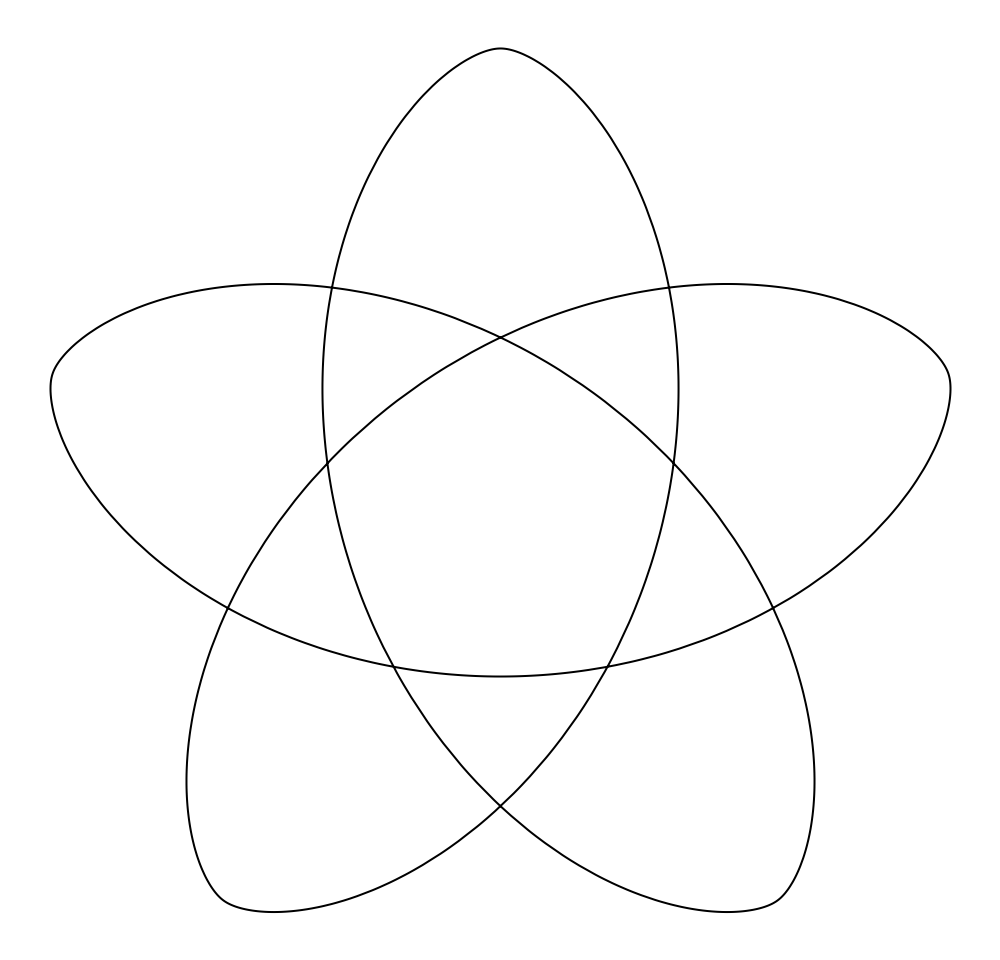}
        \caption{A $3$-rosette}
        \label{fig:Fig_rosettes02}
    \end{subfigure}
    \hfill
    \begin{subfigure}[h]{0.32\textwidth}
        \centering
        \includegraphics[width=\textwidth]{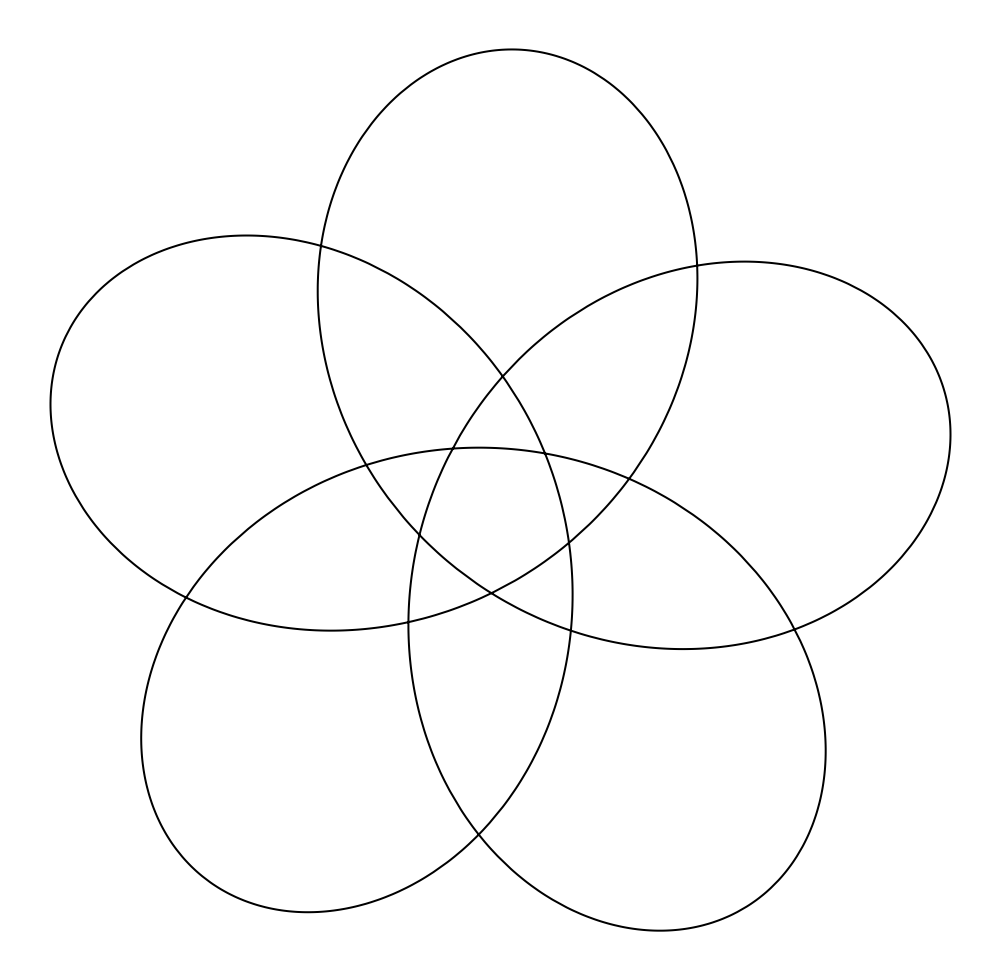}
        \caption{A $4$-rosette}
        \label{fig:Fig_rosettes03}
    \end{subfigure}
    \caption{Rosettes}
    \label{fig:Fig_rosettes}
\end{figure}

Let $\mathcal{C}$ be a \textit{smooth planar curve}, that is the image of the $C^{\infty}$ map from an interval to $\mathbb{R}^2$. A curve is \textit{regular} if its velocity vector does not vanish and is \textit{closed} if it is a $C^{\infty}$ map from $S^1$ to $\mathbb{R}^2$. A regular closed curve is \textit{locally convex} if its curvature does not vanish. A locally convex curve with a rotation number equal to $m$ is called an \text{$m$-rosette} (see Figure \ref{fig:Fig_rosettes}). An \textit{oval} is a $1$-rosette.

\begin{figure}[h]
    \centering
    \includegraphics[width=0.52137\linewidth]{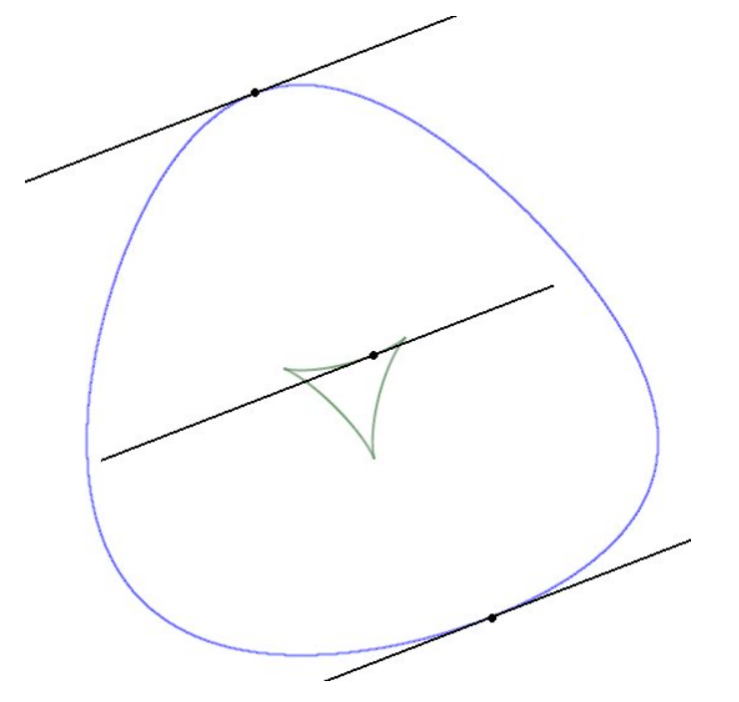}
    \caption{A curve $M$, its Wigner caustic $\E_{0.5}(M)$, and tangent lines to a parallel pair $a,b\in M$ and to $\frac{a+b}{2}\in\E_{0.5}(M)$}
    \label{fig:enter-label22}
\end{figure}

\begin{defn}
A pair $a,b\in \mathcal{C}$ ($a\neq b$) is a \textit{parallel pair} of a regular curve $\mathcal{C}$ if tangent lines to $\mathcal{C}$ at points $a,b$ are parallel (see Figure \ref{fig:enter-label22}). A chord/line through $a$ and $b$ is called an \textit{affine chord/line}.
\end{defn}

Now we will describe the most important definitions which will be relevant to future work.

\begin{defn}
Let us fix $\lambda\in\mathbb{R}$. An \textit{affine $\lambda$-equidistant} of a regular curve $\mathcal{C}$, $\Eq_\lambda(\mathcal{C})$, is the locus of points which divide an affine segment in a fixed ratio $\lambda$, that is, the following set
$$\Eq_\lambda(\mathcal{C}):=\left\{\lambda a+(1-\lambda)b\,\big|\,a,b\text{ is a parallel pair of }\mathcal{C}\right\}.$$
The \textit{Wigner caustic} of $\mathcal{C}$ is an affine $0.5$-equidistant of $\mathcal{C}$.
\end{defn}

\begin{figure}[h]
    \centering
    \begin{subfigure}[h]{0.44\textwidth}
        \centering
        \includegraphics[width=\textwidth]{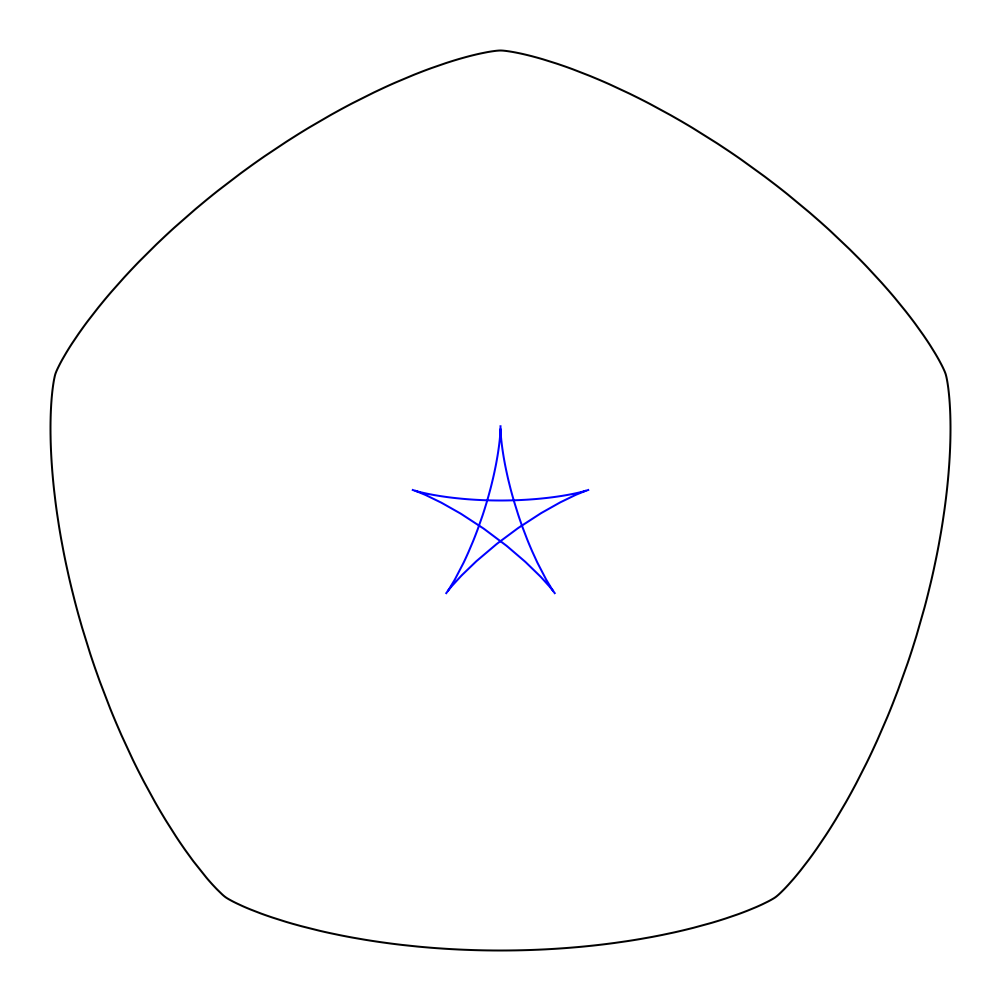}
        \caption{$\Eq_{0.5}(\mathcal{C})$}
        \label{fig:wcEqCss01}
    \end{subfigure}
    \hfill
    \begin{subfigure}[h]{0.44\textwidth}
        \centering
        \includegraphics[width=\textwidth]{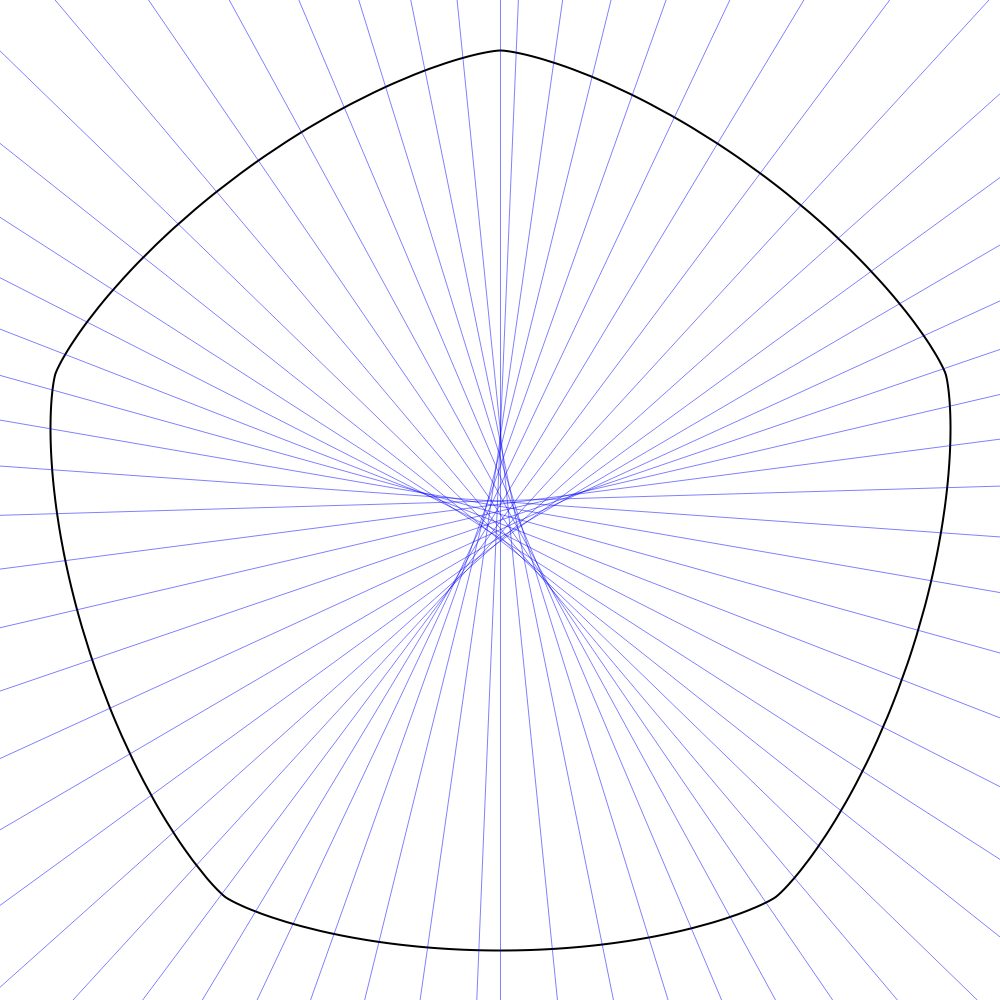}
        \caption{$\Eq_{0.5}(\mathcal{C})$ as an envelope of lines}
        \label{fig:wcEqCss02}
    \end{subfigure}
    \\ 
    \begin{subfigure}[h]{0.44\textwidth}
        \centering
        \includegraphics[width=\textwidth]{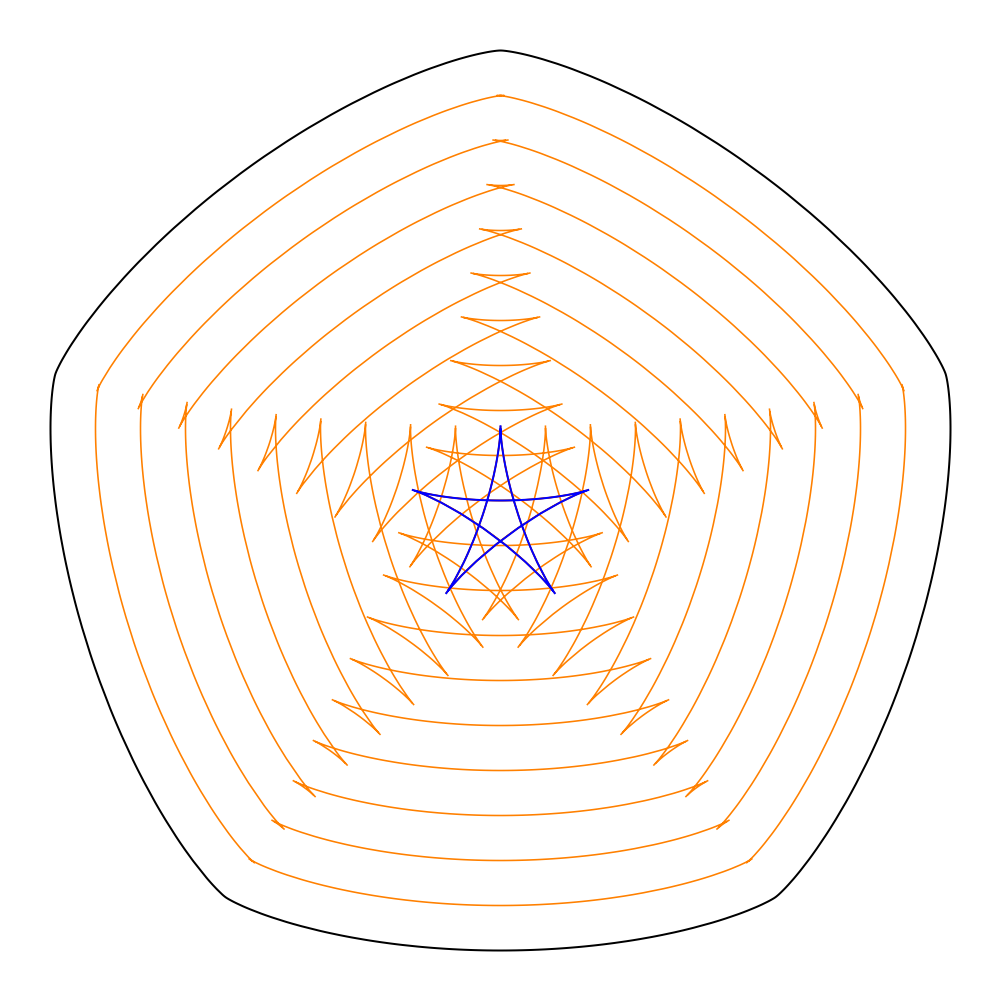}
        \caption{Affine equidistant sets of $\mathcal{C}$}
        \label{fig:wcEqCss03}
    \end{subfigure}
    \hfill
    \begin{subfigure}[h]{0.44\textwidth}
        \centering
        \includegraphics[width=\textwidth]{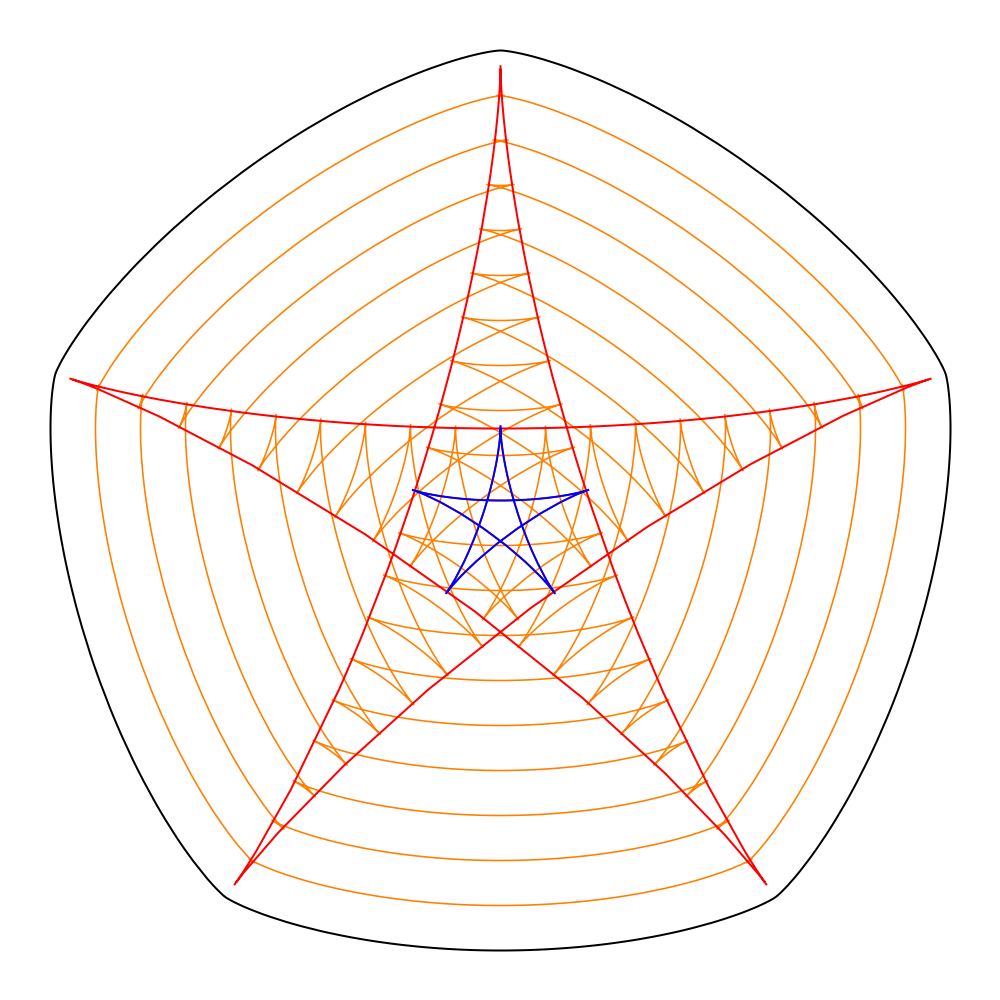}
        \caption{Affine equidistant sets of $\mathcal{C}$, $\Css(\mathcal{C})$}
        \label{fig:wcEqCss04}
    \end{subfigure}
    \hfill
    \begin{subfigure}[h]{0.44\textwidth}
        \centering
        \includegraphics[width=\textwidth]{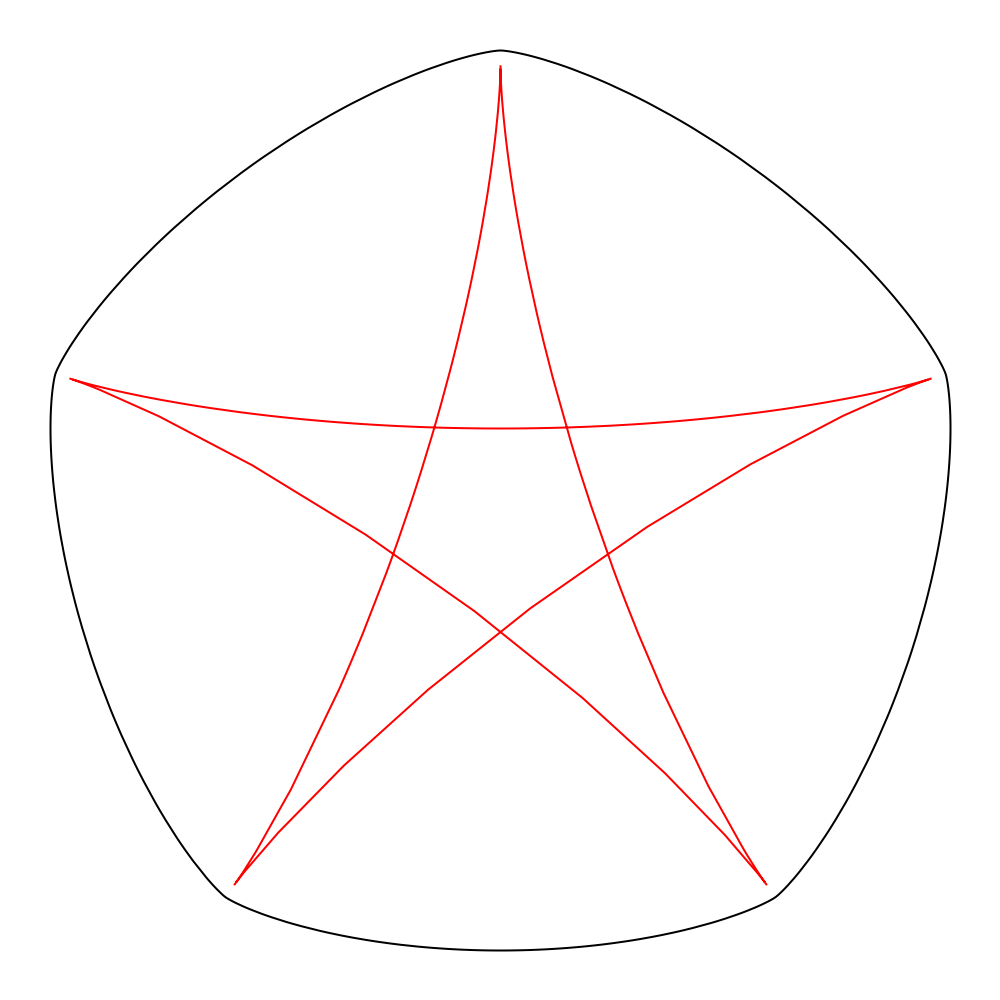}
        \caption{$\Css(\mathcal{C})$}
        \label{fig:wcEqCss05}
    \end{subfigure}
    \hfill
    \begin{subfigure}[h]{0.44\textwidth}
        \centering
        \includegraphics[width=\textwidth]{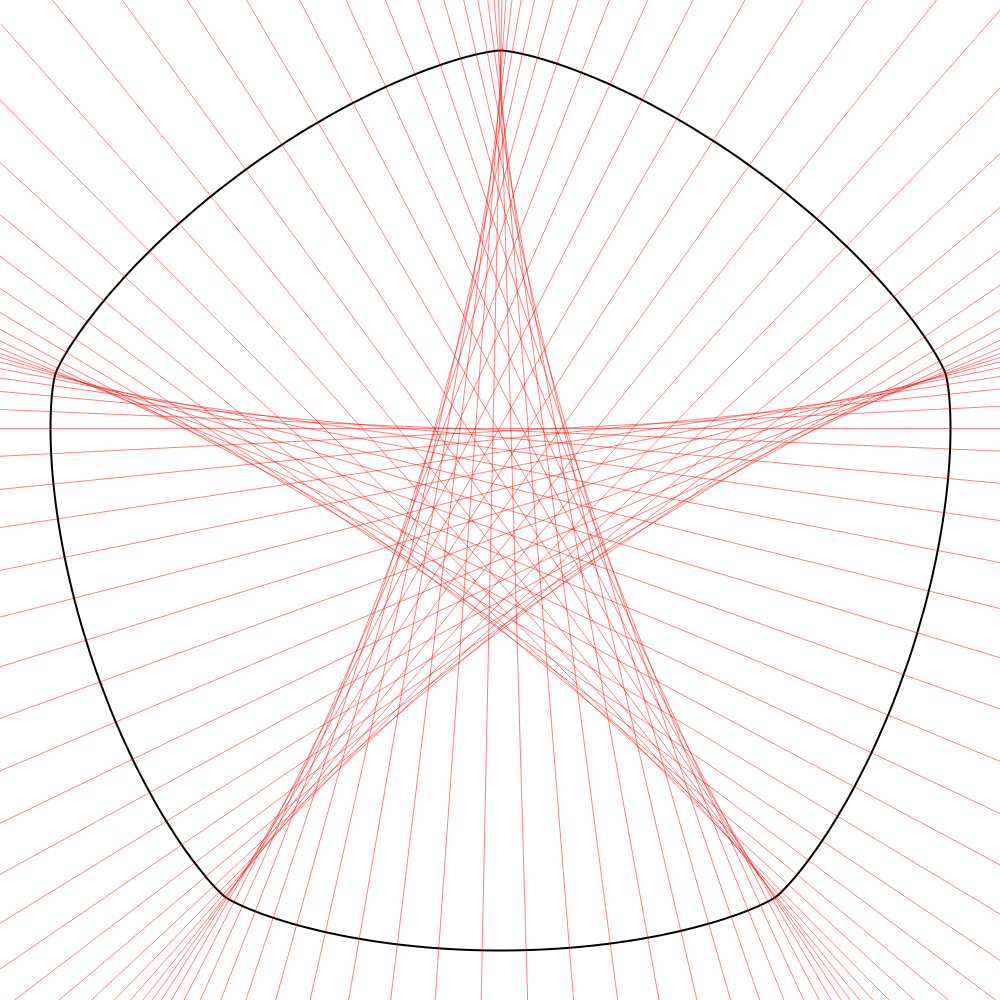}
        \caption{$\Css(\mathcal{C})$ as an envelope of affine chords}
        \label{fig:wcEqCss06}
    \end{subfigure}
    \caption{Affine sets of an oval $\mathcal{C}$}
    \label{fig:wcEqCss}
\end{figure}

If $\mathcal{C}$ is a closed regular curve, then $\Eq_{\lambda}(\mathcal{C})=\Eq_{1-\lambda}(\mathcal{C})$. Therefore, the case when $\lambda=0.5$ is in some sense special. We illustrate an oval $\mathcal{C}$ together with few affine $\lambda$-equidistants, including the Wigner caustic of $\mathcal{C}$, in Figure \ref{fig:wcEqCss}.

\begin{defn}\label{def:css_envelope}
The Centre Symmetry Set of a regular curve $\mathcal{C}$, $\Css(\mathcal{C})$, is the envelope of affine chords of $\mathcal{C}$.
\end{defn}

As we can see in Figure \ref{fig:wcEqCss}, the Centre Symmetry Set of $\mathcal{C}$ is also the set of all singular points of all affine $\lambda$-equidistants. For a generic smooth regular planar curve $\mathcal{C}$, affine $\lambda$-equidistant (for a generic value of $\lambda$) and the Centre Symmetry Set of $\mathcal{C}$ consists of smooth curves (\textit{branches}) with at most cusp singularities. Each such branch of $\Eq_{\lambda}(\mathcal{C})$ has even number of cusps provided that $\lambda\neq 0.5$ (\cite{DZ-WC}). If $\mathcal{C}$ is an oval, the number of cusps of the Wigner caustic and the number of cusps of the Centre Symmetry Set are odd (\cite{B1, GH1}), and the number of cusps of $\Eq_{0.5}(\mathcal{C})$ is not smaller than the number of cusps of $\Css(\mathcal{C})$ (\cite{DR1}). See \cite{DZ-singular, DZ-WC, Schneider1, Schneider2, Zwierz1, Zwierz2, Zwierz3} for theorems on global geometry of the Wigner caustic and other affine equidistants of a planar curve, and see \cite{DZGauss, GH1} for theorems on global geometry of the Centre Symmetry Set of a planar curve. For details in higher dimensions see \cite{DJRR1, DMR1, DR1} and the literature therein.

In differential geometry and singularity theory, mathematicians are looking for global invariants, and global theorems concerning of sets. In \cite{DZ-singular} Authors proved a global property for the Wigner caustic of a loop.

\begin{defn}\label{def:eq_envelope}
A simple smooth curve $\gamma:(s_1,s_2)\to\mathbb{R}^2$ with non-vanishing curvature is called a \textit{loop} if $\displaystyle\lim_{s\to s_1^+}\gamma(s)=\lim_{s\to s_1^-}\gamma(s)$. A loop $\gamma$ is called \textit{convex} if the absolute value of its rotation number is not greater than $1$, otherwise it is called \textit{non-convex}.
\end{defn}

\begin{figure}[h]
    \centering
    \includegraphics[scale=0.4]{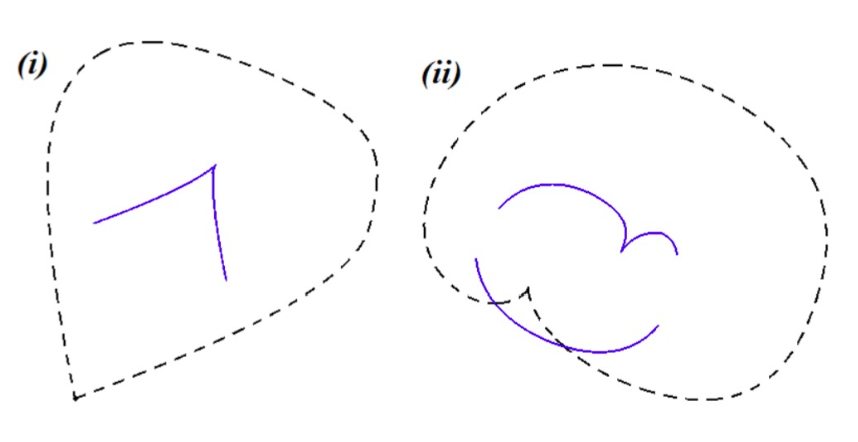}
    \caption{(i) a convex loop and its Wigner caustic, (ii) a non-convex loop and its Wigner caustic}
    \label{fig:loops}
\end{figure}

We illustrate examples of loops in Figure \ref{fig:loops}.

\begin{thm}[\cite{DZ-singular}]\label{TheoremLoops} The Wigner caustic of a loop has a singular point.
\end{thm}

\begin{figure}[h]
    \centering
    \includegraphics[scale=0.4]{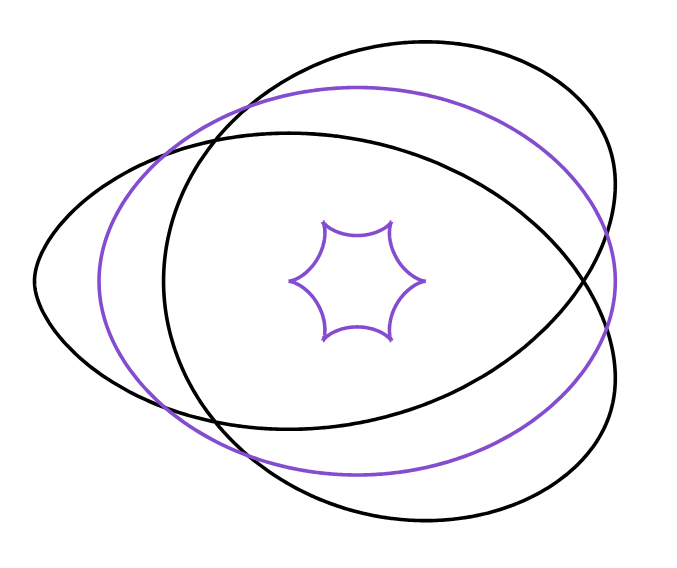}
    \caption{A $2$-rosette and its Wigner caustic}
    \label{fig:2rosetteconvexloops}
\end{figure}

\begin{figure}[h]
    \centering
    \includegraphics[scale=0.4]{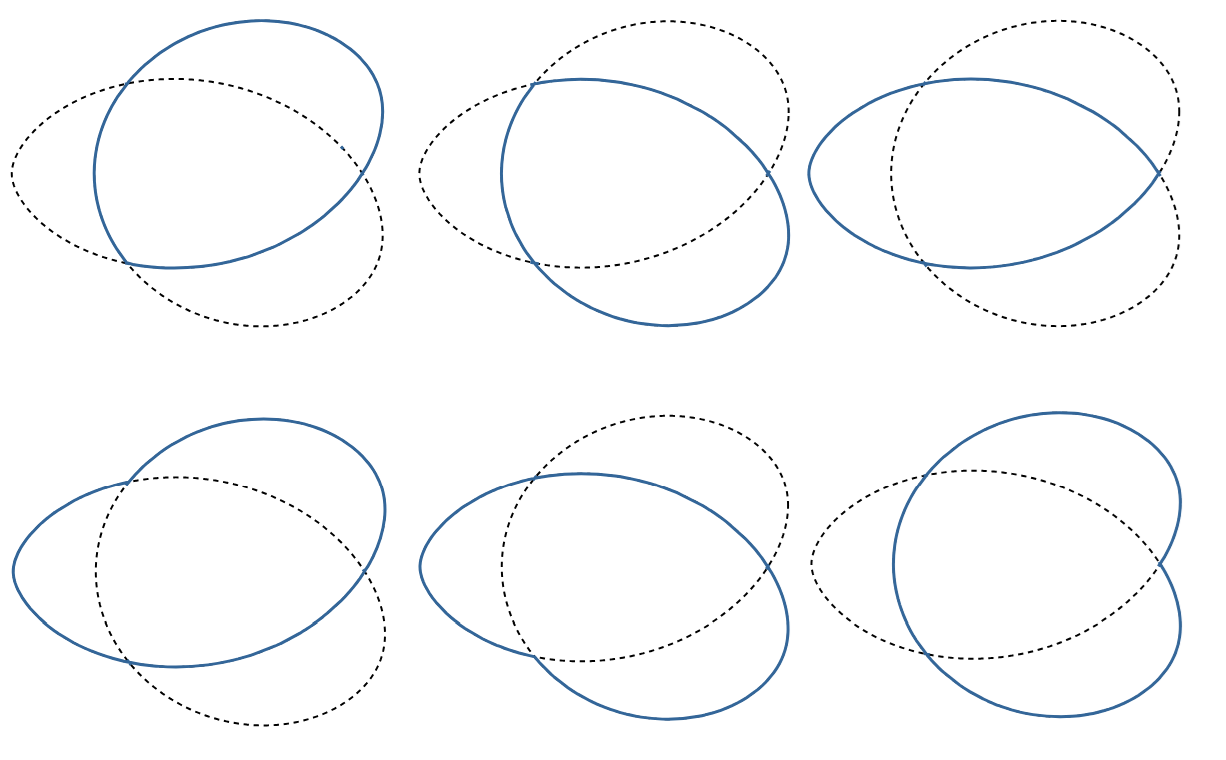}
    \caption{Convex and non-convex loops of curve presented in Figure \ref{fig:2rosetteconvexloops}}
    \label{fig:enter-label}
\end{figure}

Now we will illustrate an application of Theorem \ref{TheoremLoops}. In Figure \ref{fig:2rosetteconvexloops} we present a $2$-rosette and its Wigner caustic. Since this curve has exactly three convex and three non-convex loops (as shown in Figure \ref{fig:enter-label}), its Wigner caustic has at least $6$ singular points. In this case, Theorem \ref{TheoremLoops} helps us localize all singular points of the Wigner caustic. Also, see Figure \ref{fig: hearts} (Hearts), where we can see singular points of the Wigner caustic which come from the loops.

Now we will present a proposition that will help us define an affine $\lambda$-equidistant as an envelope of lines. Let $T_xX$ denote the tangent space to $X$ at $x\in X$.

\begin{prop}[\cite{DZ-WC}]\label{prop-envelope} Let $a,b$ be a parallel pair of a regular closed curve $\mathcal{C}$. If $p:=\lambda a+(1-\lambda)b$ is a regular point of $\Eq_\lambda(\mathcal{C})$. Then lines $T_a\mathcal{C}$, $T_b\mathcal{C}$, and $T_p\Eq_\lambda(\mathcal{C})$ are parallel.
\end{prop}

See Figure \ref{fig:enter-label22} as an illustration of Proposition \ref{prop-envelope}.

\begin{col}\label{CorEqEnvelope}
Let $\mathcal{C}$ be a regular closed curve. Then an affine $\lambda$-equidistant is the envelope of the following family of lines
\begin{align*}
    \left\{\lambda T_aM+(1-\lambda)T_bM\ \big|\ a,b\text{ is a parallel pair of }\mathcal{C}\right\},
\end{align*}
where $+$ is the Minkowski addition, and $aX:=\{ax\,|\,x\in X\}$.
\end{col}

We illustrate Corollary \ref{CorEqEnvelope} in Figure \ref{fig:wcEqCss02}.

In general, drawing the Centre Symmetry Set and affine equidistants requires finding parallel pairs of a regular curve $\mathcal{C}$. It can be done by solving a large number of differential equations, and then plotting appropriate lines or points that will form the sets under consideration. This approach is burdened with approximation errors of algorithms solving differential equations, inconsistency of these methods, and high computational complexity. 

As in \cite{Zwierz1}, we will use a different method to find the parallel pair of a curve. We will use the so-called \textit{Minkowski support function}. This method works for the \textit{hedgehogs} (see \cite{MM1,MM2} and the literature therein), i.e. curves that can be parameterized using their Gauss maps. Regular closed hedgehogs are rosettes.

Let $R_n$ be an $n$-rosette and let $\boldsymbol 0$ be the origin of the plane $\mathbb{R}^2$. Let
$$[0,2n\pi)\ni\theta\mapsto\gamma(\theta)\in\mathbb{R}^2$$
be a parameterization of $R_n$ in terms of the tangential angle $\theta$ to the rosette $R_n$. Let $[0,2n\pi)\ni\theta\mapsto p(\theta)\in\mathbb{R}$ be a function, where $p(\theta)$ for a fixed value of $\theta$ is an oriented distance between $\boldsymbol 0$ and the tangent line to $R_n$ at a point $\gamma(\theta)$ in the direction $(\cos\theta,\sin\theta)$. In terms of the \textit{polar-tangential coordinates} we can find that the parameterization $\gamma$ has the following form:
\begin{align}
\label{eq:support_fun}[0,2n\pi)\ni\theta\mapsto\gamma(\theta)=\big(p(\theta)\cos\theta -p'(\theta)\sin\theta, p(\theta)\sin\theta + p'(\theta)\cos\theta\big)\in\mathbb{R}^2.
\end{align}
This parameterization is very useful for finding the parallel pairs -- notice that the pair $\gamma(\theta),\gamma(\theta+k\pi)$ is a parallel pair for $k=1, 2, \ldots, n-1$. Furthermore, the curvature function is given by the formula
$$\kappa(\theta)=\dfrac{1}{p(\theta)+p''(\theta)}>0.$$
From \cite{DZGauss, DZ-WC, Zwierz3} we get that the Wigner caustic of $R_n$ consists of exactly $n$ smooth branches:
\begin{itemize}
    \item a branch $\Eq_{0.5, k}(R_n)$ for $k=1,\ldots, n-1$ with the parameterization
    \begin{align}
    \label{eq:wcparam_1}    [0,2n\pi)\ni\theta\mapsto\dfrac{1}{2}\big(\gamma(\theta)+\gamma(\theta+k\pi)\big)\in\mathbb{R}^2,
    \end{align}
    \item a branch $\Eq_{0.5,n}(R_n)$ with the parameterization
    \begin{align}
    \label{eq:wcparam_1}  [0,n\pi)\ni\theta\mapsto\dfrac{1}{2}\big(\gamma(\theta)+\gamma(\theta+n\pi)\big)\in\mathbb{R}^2,
    \end{align}
\end{itemize}
the affine $\lambda$-equidistant of $R_n$, where $\lambda\neq 0,\frac{1}{2},1$, consists of exactly $2n-1$ smooth branches:
\begin{itemize}
    \item a branch $\Eq_{\lambda,k}(R_n)$ for $k=1,\ldots,n-1$ with the parameterization
    \begin{align} \label{eq:eqq_1}
    [0,2n\pi)\ni\theta\mapsto\lambda\gamma(\theta)+(1-\lambda)\gamma(\theta+k\pi)\in\mathbb{R}^2,
    \end{align}
    \item a branch $\Eq_{\lambda,n}(R_n)$ with the parameterization
    \begin{align} \label{eq:eqq_2}
    [0,2n\pi)\ni\theta\mapsto\lambda\gamma(\theta)+(1-\lambda)\gamma(\theta+n\pi)\in\mathbb{R}^2,
    \end{align}
    \item a branch $\Eq_{\lambda,n+k}(R_n)$ for $k=1,\ldots,n-1$ with the parameterization
    \begin{align} \label{eq:eqq_3}
    [0,2n\pi)\ni\theta\mapsto(1-\lambda)\gamma(\theta)+\lambda\gamma(\theta+k\pi)\in\mathbb{R}^2,
    \end{align}
\end{itemize}
the Centre Symmetry Set of $R_n$ consists of $n$ branches:
\begin{itemize}
    \item a branch $\Css_{k}(R_n)$ for $k=1,\ldots, n-1$ with the parameterization
    \begin{align}\label{eq:css_1}
    [0,2n\pi)\ni\theta\mapsto\dfrac{\kappa(\theta)\gamma(\theta)+(-1)^k\kappa(\theta+k\pi)\gamma(\theta+k\pi)}{\kappa(\theta)+(-1)^k\kappa(\theta+k\pi)}\in\mathbb{R}^2,
    \end{align}
    \item a branch $\Css_{n}(R_n)$ with the parameterization
    \begin{align}\label{eq:css_2}
    [0,n\pi)\ni\theta\mapsto\dfrac{\kappa(\theta)\gamma(\theta)+(-1)^n\kappa(\theta+n\pi)\gamma(\theta+n\pi)}{\kappa(\theta)+(-1)^n\kappa(\theta+n\pi)}\in\mathbb{R}^2.
    \end{align}
\end{itemize}

\begin{figure}[h]
    \centering
    \begin{subfigure}[h]{0.44\textwidth}
        \centering
        \includegraphics[width=\textwidth]{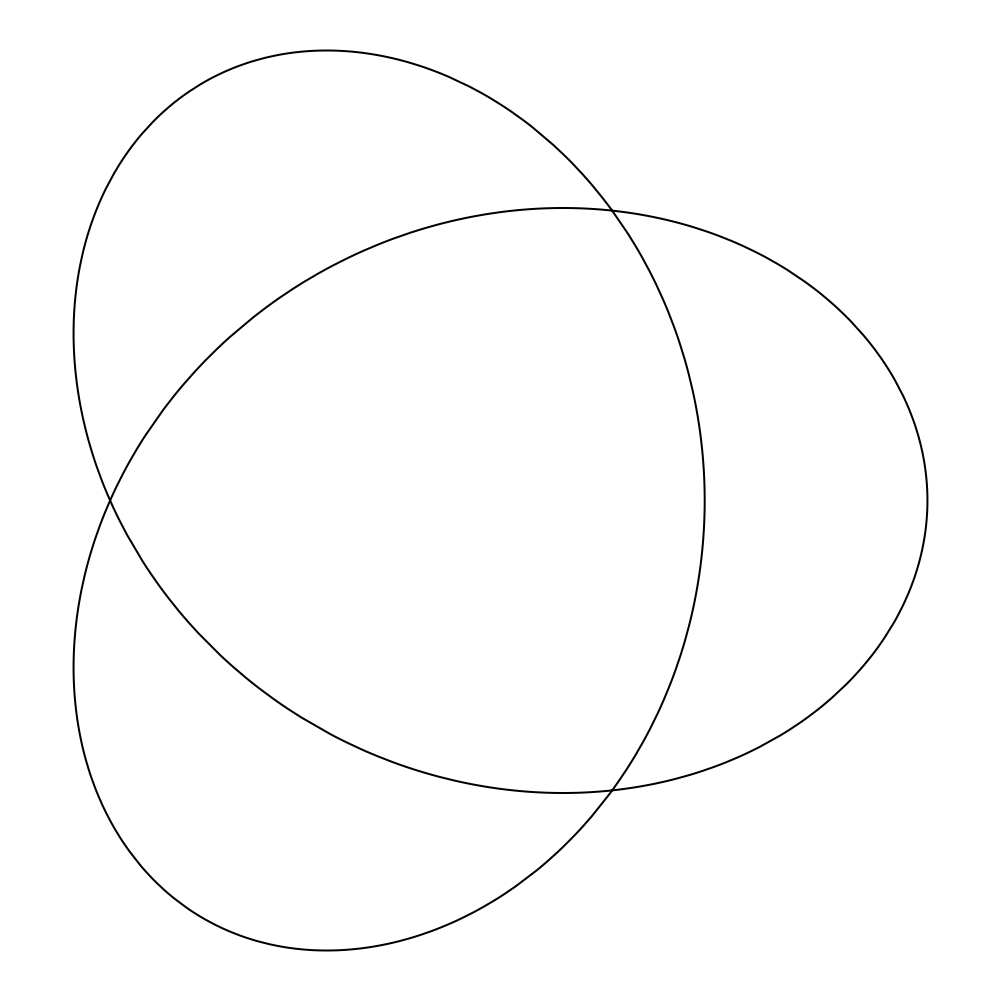}
        \caption{A $2$-rosette $\mathcal{R}$}
        \label{fig:Fig_Rosette_wceqcss01}
    \end{subfigure}
    \hfill
    \begin{subfigure}[h]{0.44\textwidth}
        \centering
        \includegraphics[width=\textwidth]{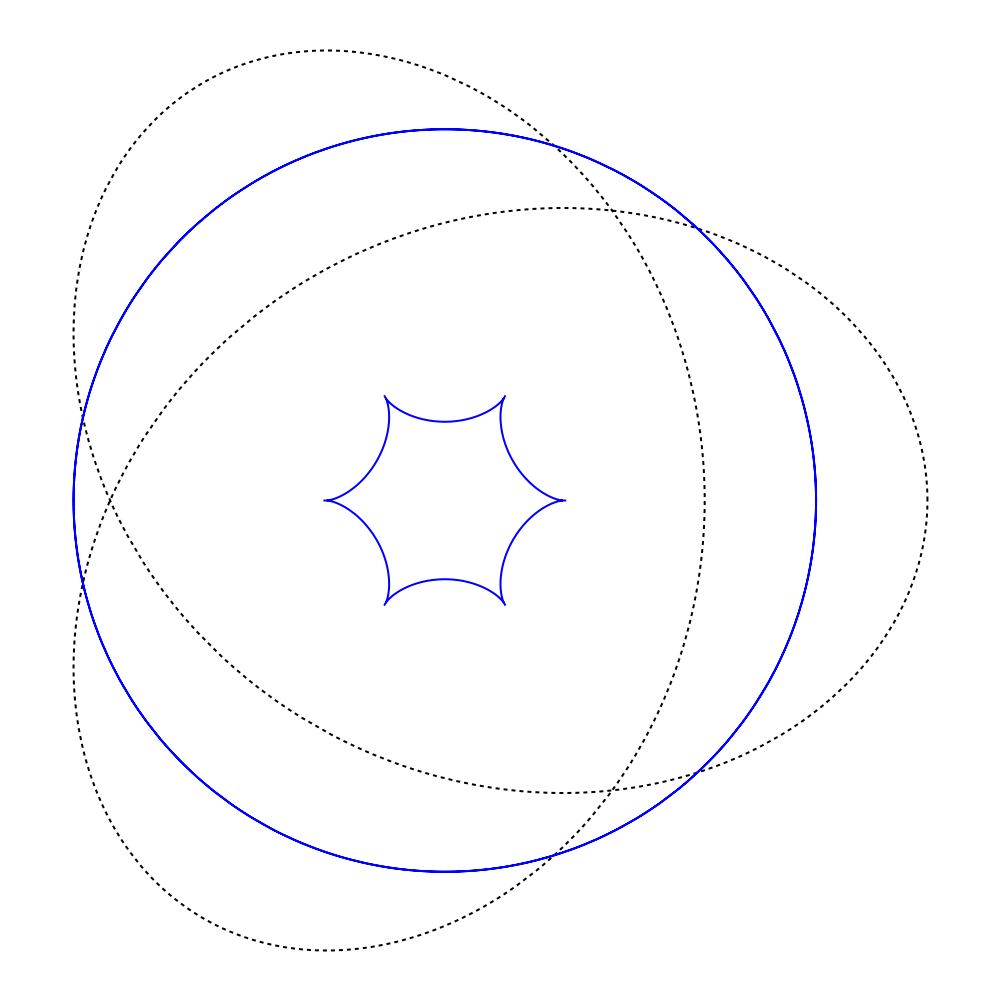}
        \caption{$\Eq_{0.5}(\mathcal{R})$}
        \label{fig:Fig_Rosette_wceqcss02}
    \end{subfigure}
    \\ 
    \begin{subfigure}[h]{0.44\textwidth}
        \centering
        \includegraphics[width=\textwidth]{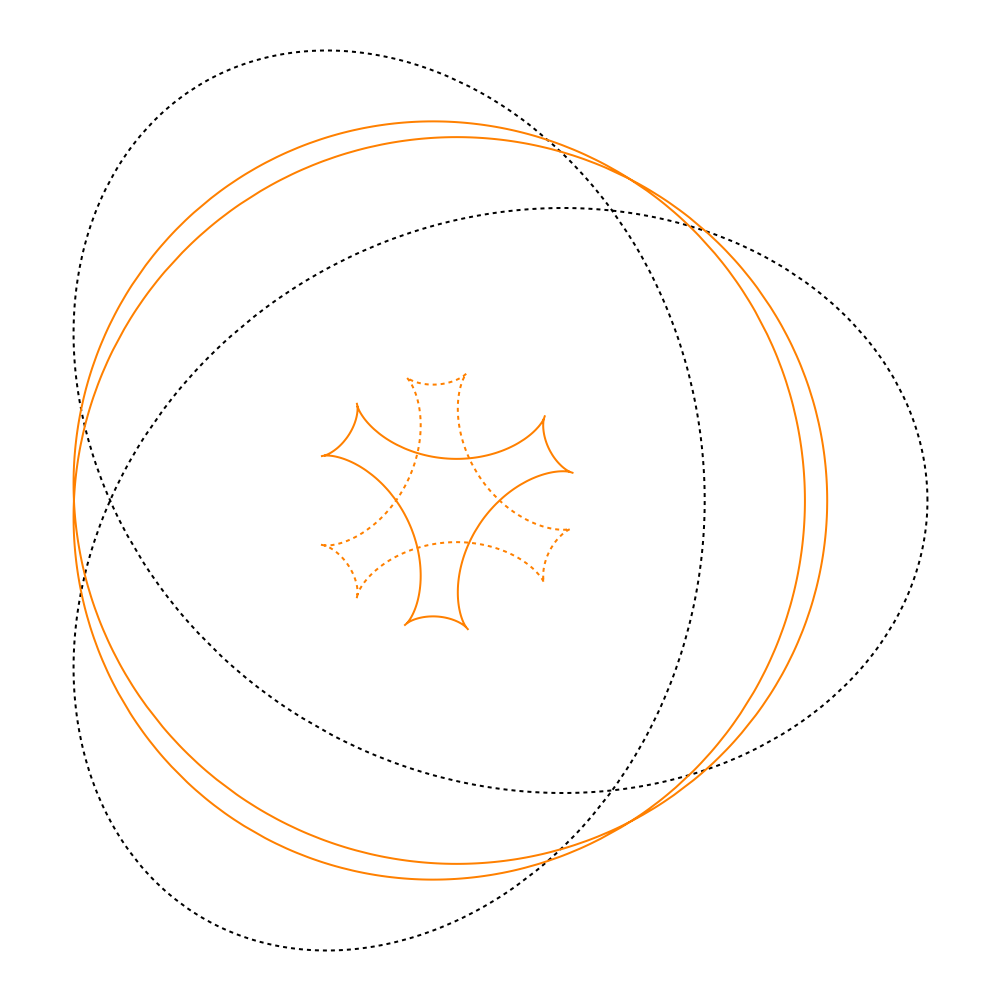}
        \caption{$\Eq_{0.45}(\mathcal{R})$}
        \label{fig:Fig_Rosette_wceqcss03}
    \end{subfigure}
    \hfill
    \begin{subfigure}[h]{0.44\textwidth}
        \centering
        \includegraphics[width=\textwidth]{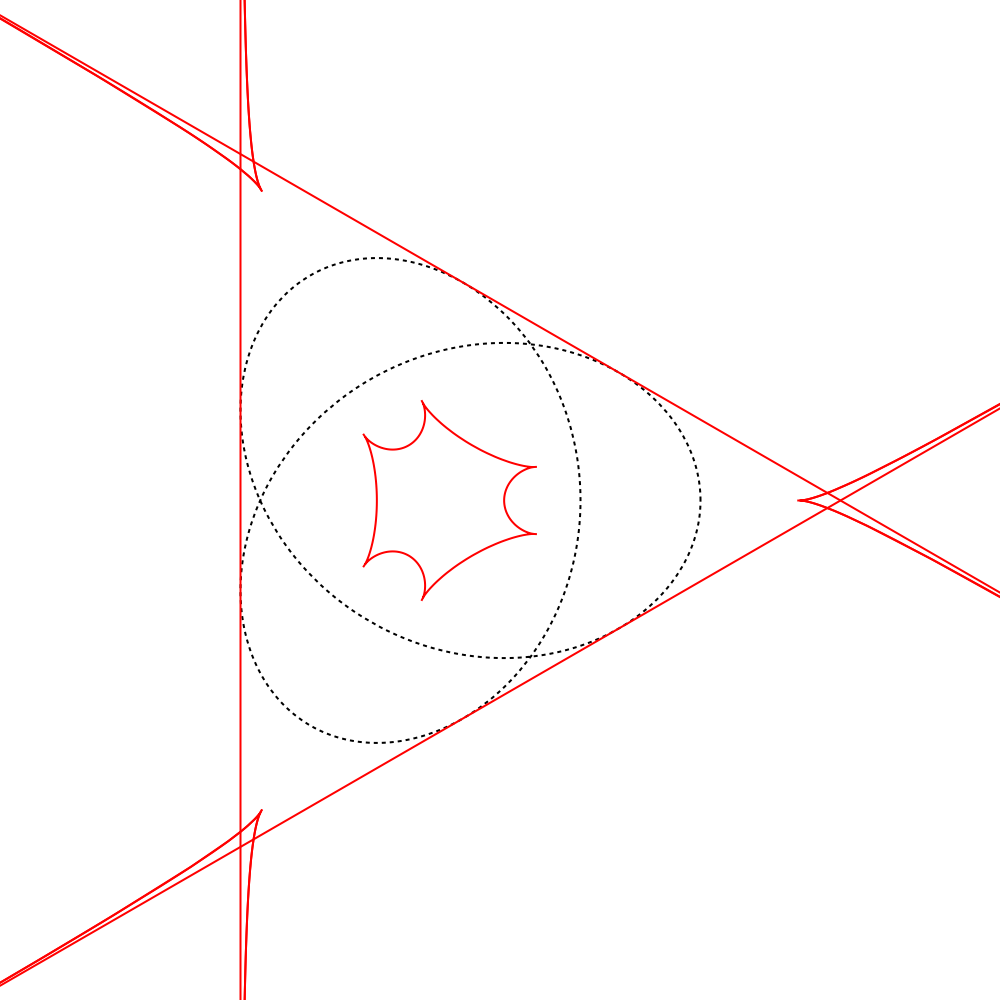}
        \caption{$\Css(\mathcal{R})$}
        \label{fig:Fig_Rosette_wceqcss04}
    \end{subfigure}
    \caption{Affine sets of a rosette}
    \label{fig:Fig_Rosette_wceqcss}
\end{figure}

For details and the proofs of other geometric properties of the branches above see \cite{DZGauss, DZ-WC, Zwierz3}. In Figure \ref{fig:Fig_Rosette_wceqcss} we illustrate the branches of the Wigner caustic, the branches of the affine $0.45$-equidistant, and the branches of the Centre Symmetry Set of the $2$-rosette $\mathcal{R}$ with the support function $p(\theta)=10+3\cos\frac{3t}{2}$. Note that the bi-tangent lines to $\mathcal{R}$ are subset of $\Css(\mathcal{R})$ (\cite{GH1}).

\section{Lapidary Construction of Image in \texttt{Mathematica}}

In this section, we present a simplified construction of figures containing the Center Symmetry Set, $\lambda$-equidistant sets, including Wigner caustic, using Wolfram Mathematica (\cite{WMath}).

For a fixed support function $[0,2n\pi)\ni\theta\mapsto p(\theta)\in\mathbb{R}$ of an $n$-rosette $R_n$ there exists a useful parametrization of $R_n$, the branches of affine equidistants (including the Wigner caustic), and the Centre Symmetry Set of $R_n$ given by the equations from \eqref{eq:support_fun} to \eqref{eq:css_2}. Furthermore, we can adopt these equations for the union of two, or more, rosettes to find the parameterizations of branches of the appropriate sets. By Definitions \ref{def:css_envelope} and \ref{def:eq_envelope} we can illustrate these sets as the envelopes of the family of lines.  

By setting the appropriate curves or a single curve, thickness, transparency, colors of a family of lines, background, and many other seemingly insignificant parameters, we obtained many interesting artistic visualizations of the sets.

Now, let us analyze one example of a lapidary construction of an artistic visualization of studied sets, using Wolfram Mathematica.

We start by defining the support function (line 1 in Listing \ref{list1}). Then, using the transformation of the support points, we determine the parametrization of the rosette $\gamma$ (line 2 in Listing \ref{list1}).
\begin{lstlisting}[caption=Definition od support function and curve gamma, label = list1]
supportfunction[t_] = 150 + 4*Sin[5*t] + Cos[4*t]
gamma[t_] = {supportfunction[t]*Cos[t] - supportfunction'[t]*Sin[t], 
   supportfunction[t]*Sin[t] + supportfunction'[t]*Cos[t]};

\end{lstlisting}

Directly by the Definition of the support function one can easily get that the tangent line to a rosette and any affine equidistant of the rosette at the point corresponding to the angle $ \theta$ is parallel to the vector $(-\sin\theta, \cos \theta)$ (see Listing \ref{list2}).

\begin{lstlisting}[caption=Definition of tangent, label = list2]
tang[x_] = {-Sin[x], Cos[x]}
\end{lstlisting}

In Listing \ref{list3}, we plot the family of tangents to the curve $\gamma$. 

\begin{lstlisting}[caption= Plotting the family of tangents to the curve $\gamma$, label = list3]
plotgammablack = 
 ParametricPlot[gamma[t], {t, 0, 2 Pi}, Axes -> False, 
 PlotStyle -> {Black, Thick}, Background -> Black]
\end{lstlisting}

\begin{figure}[h]
    \centering
    \includegraphics[scale=0.3]{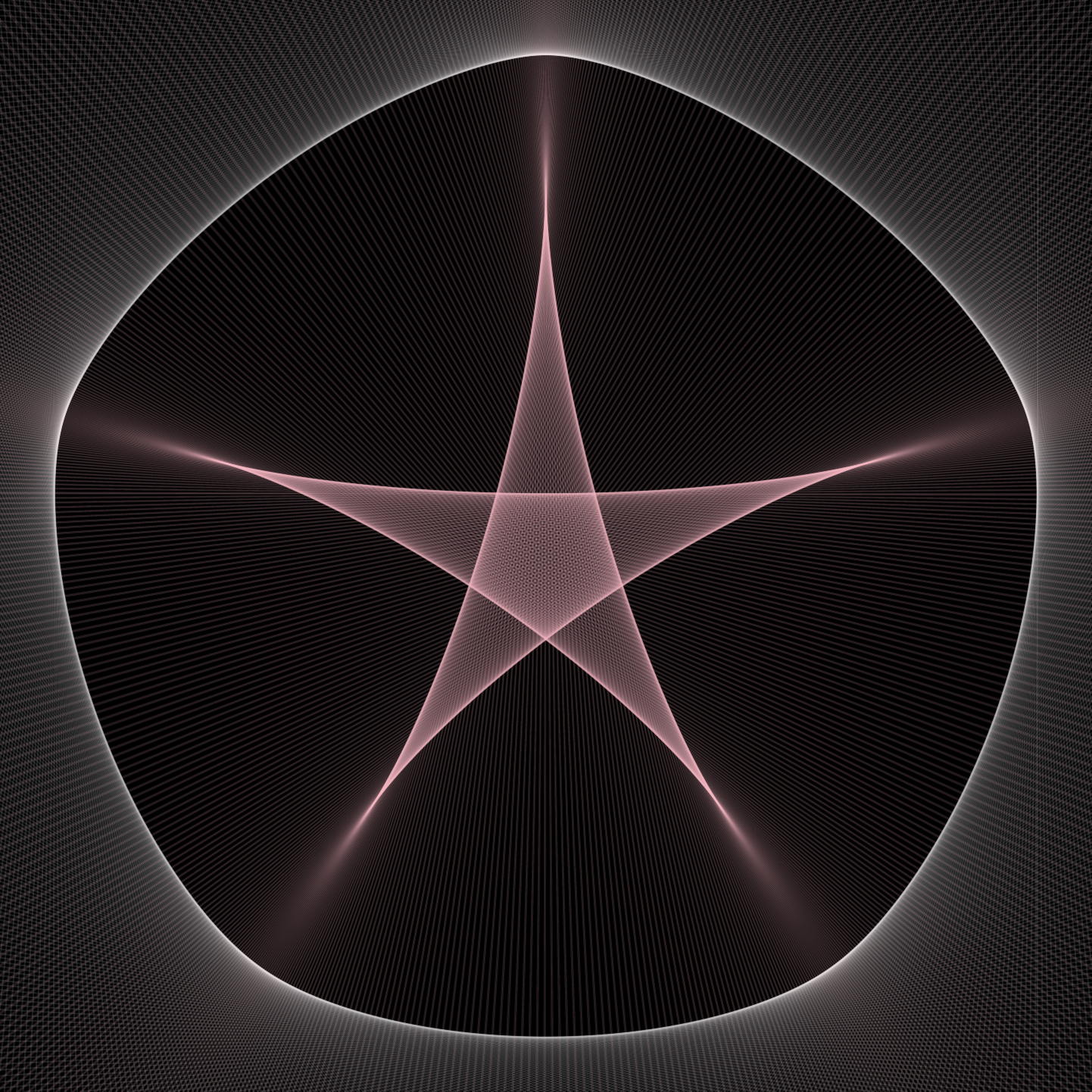}
    \caption{Graph of Center Symmetry Set (CSS)}
    \label{fig:r2}
\end{figure}

Based on the following code, we obtain the picture in  Figure \ref{fig:r2}. 
\begin{lstlisting}[caption=Ploting the family of tangent lines to $\gamma$ and the family of chords of $\gamma$ (which are tangent lines to $\Css(\gamma)$]
Show [
 plotgammablack,
 Graphics[{White, Opacity[0.1], 
   Table[InfiniteLine[gamma[s], tang[s]], {s, 0, 2 Pi, 0.01}]}],
 Graphics[{RGBColor["#FFC0CB"], Opacity[0.1], 
   Table[InfiniteLine [{gamma[a], gamma[a + Pi]}], {a, 0, 2 Pi, 
     0.01}]}]
\end{lstlisting}

In the next step, we create the family of tangent lines to few affine $\lambda$-equidistant sets for $\lambda$ values of $\frac{1}{2}$, $\frac{10}{22}, \frac{10}{24}, \frac{10}{26}$, and $ \frac{10}{28}$ -- see Listing \ref{list7}.

\begin{figure}[h]
    \centering
    \includegraphics[width=\textwidth]{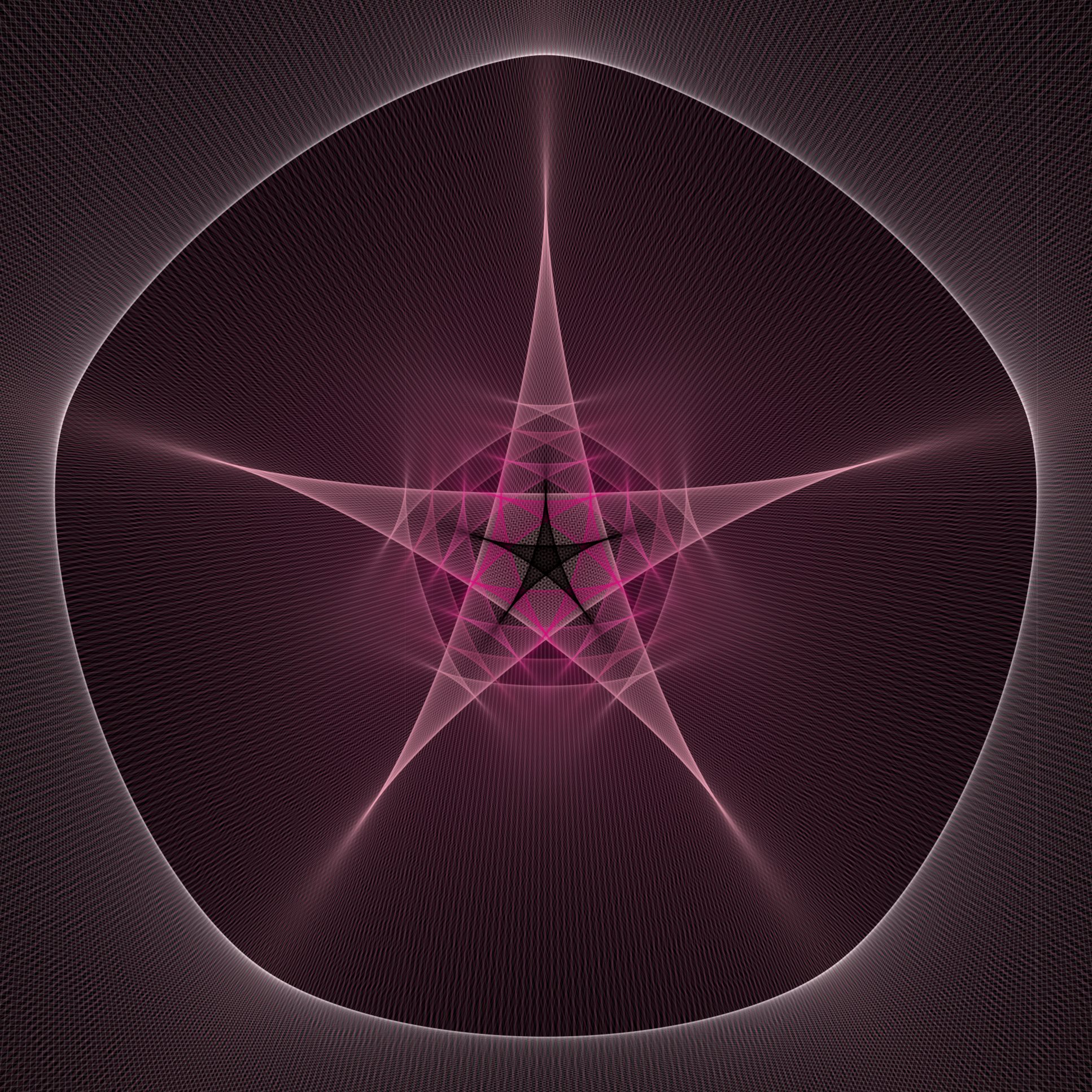}
    \caption{Final result containing CSS, the Wigner caustic and $\lambda$-equidistant sets of $\gamma$, titled \textit{Pink Star}}
    \label{fig:r3}
\end{figure}

As a result, we obtain the following outcome presented in Figure \ref{fig:r3}.

\begin{lstlisting}[caption= Drawing the family of tangent lines to few  $\lambda$-equidistant sets of $\gamma$, label = list7]
Show [
 plotgammablack,
 Graphics[{White, Opacity[0.1], 
   Line[Table[{startlinetang[s], endlinetang[s]}, {s, 0, 2 Pi, 
      0.01}]]}],
 
 Graphics[{RGBColor["#FFC0CB"], Opacity[0.1], 
   Table[InfiniteLine [{gamma[a], gamma[a + Pi]}], {a, 0, 2 Pi, 
     0.01}]}],
 Graphics[{RGBColor["#FFC0CB"], Opacity[0.05], 
   Table[InfiniteLine[affline[a, 10/28], tang[a]], {a, 0, 2 Pi, 
     0.01}]}],
 Graphics[{RGBColor["#DB7093"], Opacity[0.05], 
   Table[InfiniteLine[affline[a, 10/26], tang[a]], {a, 0, 2 Pi, 
     0.01}]}],
 Graphics[{RGBColor["#FF69B4"], Opacity[0.05], 
   Table[InfiniteLine[affline[a, 10/24], tang[a]], {a, 0, 2 Pi, 
     0.01}]}],
 Graphics[{RGBColor["#FF1493"], Opacity[0.05], 
   Table[InfiniteLine[affline[a, 10/22], tang[a]], {a, 0, 2 Pi, 
     0.01}]}],
 Graphics[{Black, Opacity[0.05], 
   Table[InfiniteLine[affline[a, 0.5], tang[a]], {a, 0, 2 Pi, 0.01}]}]
\end{lstlisting}

\section{Artistic Visualization}

\begin{figure}[h]
    \centering
    \begin{subfigure}[h]{0.44\textwidth}
        \centering
        \includegraphics[width=\textwidth]{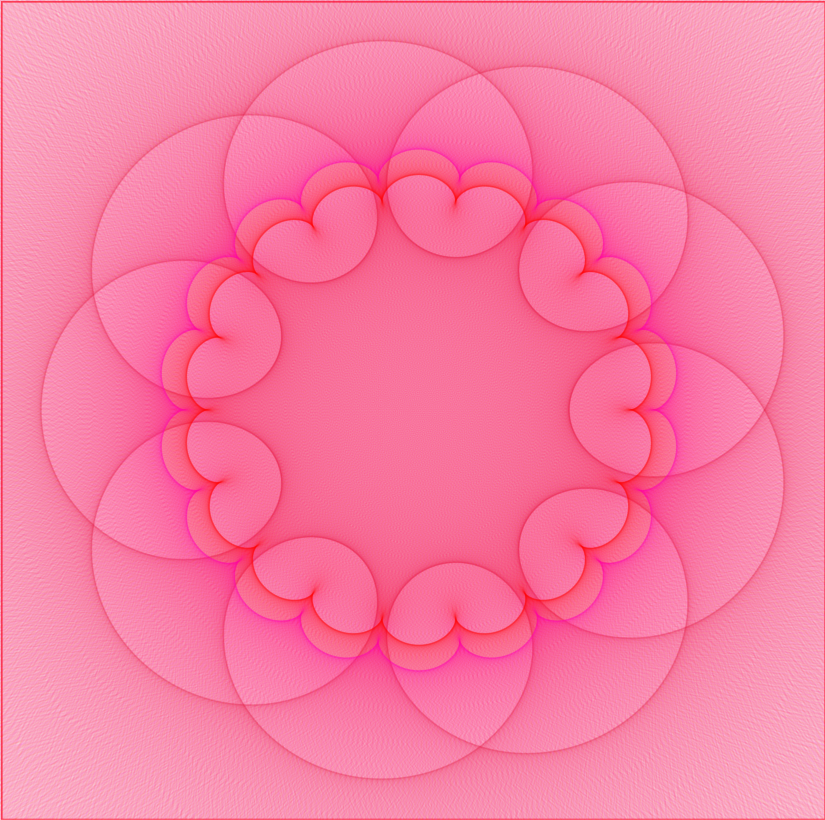}
        \caption{\textit{Hearts}}
        \label{fig: hearts}
    \end{subfigure}
    \hfill
    \begin{subfigure}[h]{0.44\textwidth}
        \centering
        \includegraphics[width=\textwidth]{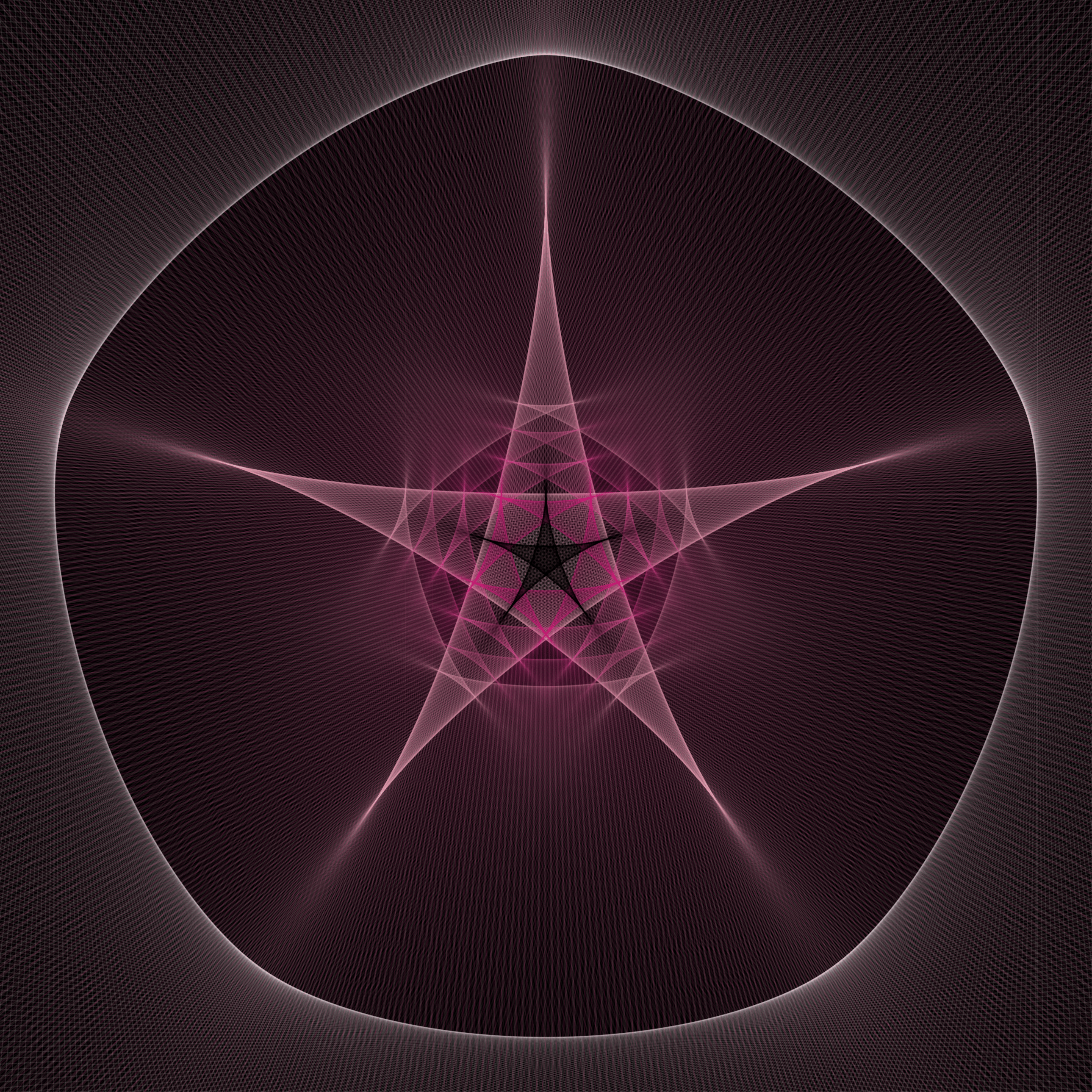}
        \caption{\textit{Pink Star}}
        \label{fig: pinkstars}
    \end{subfigure}
    \\ 
    \begin{subfigure}[h]{0.44\textwidth}
        \centering
        \includegraphics[width=\textwidth]{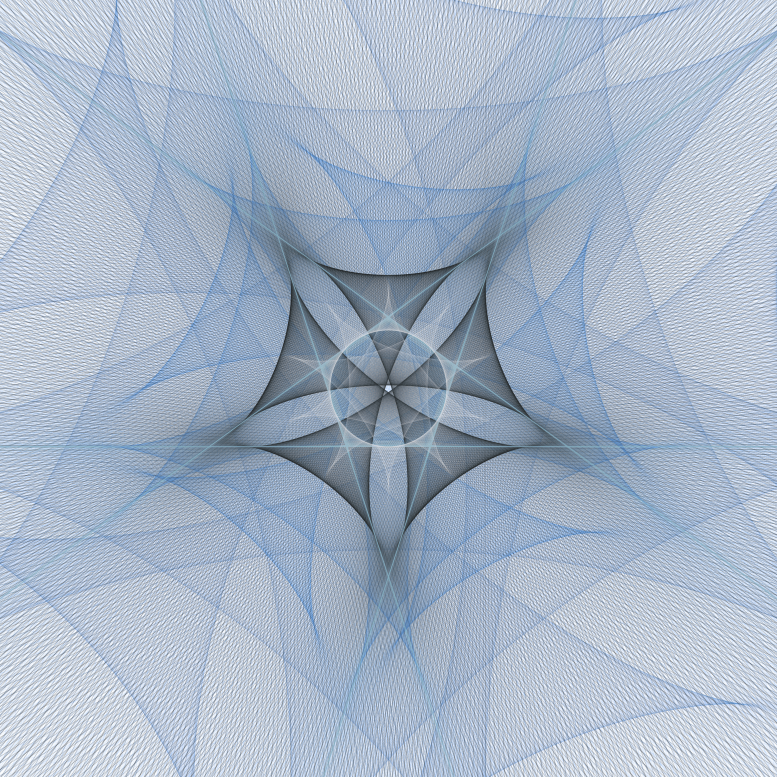}
        \caption{\textit{Small Fat Pentagram}}
        \label{fig: smallfat}
    \end{subfigure}
    \hfill
    \begin{subfigure}[h]{0.44\textwidth}
        \centering
        \includegraphics[width=\textwidth]{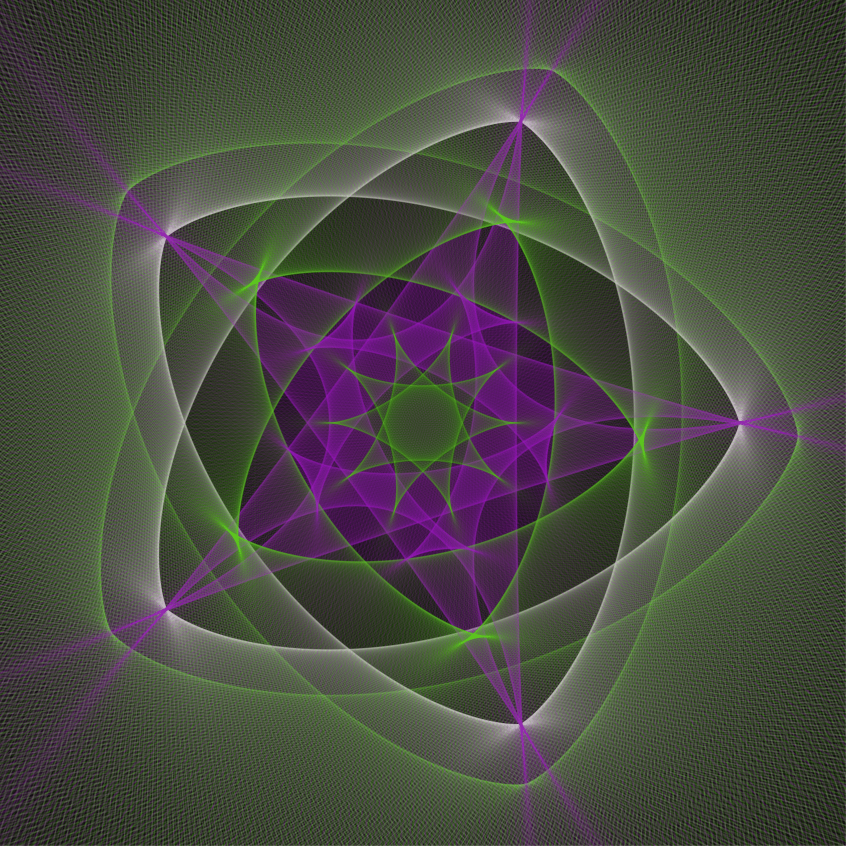}
        \caption{\textit{Big Purple Pentagramobile}}
        \label{fig: bigpurple}
    \end{subfigure}
    \label{fig:dpsz}
    \begin{subfigure}[h]{0.44\textwidth}
        \centering
        \includegraphics[width=\textwidth]{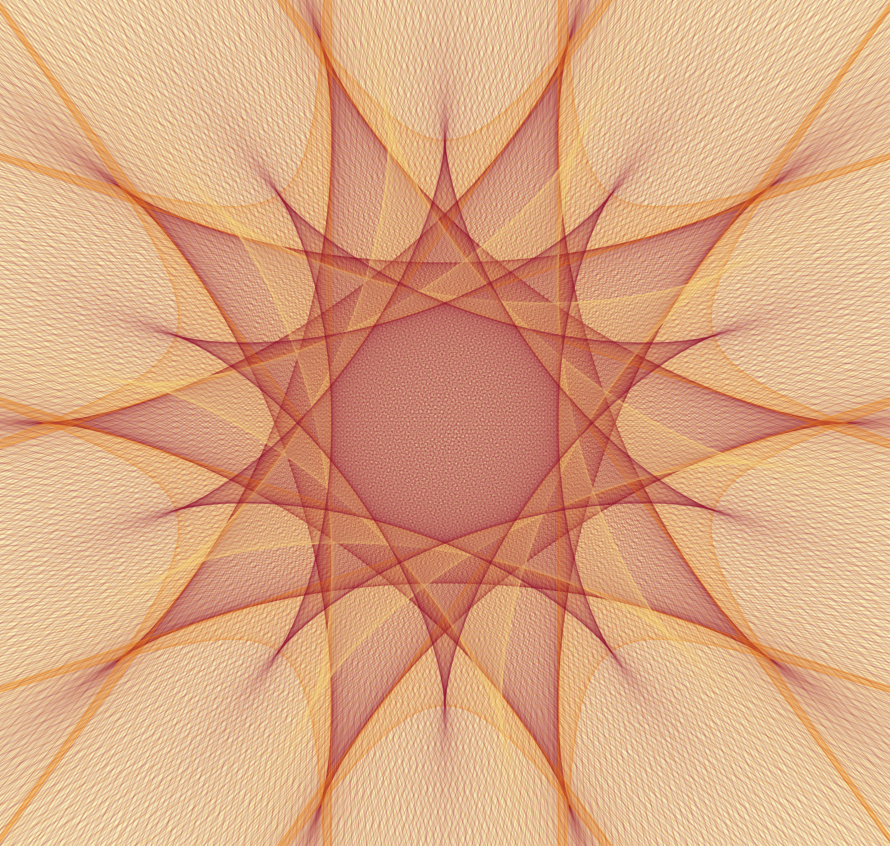}
        \caption{\textit{Charmander's Tail}}
        \label{fig: charmander}
    \end{subfigure}
    \hfill
    \begin{subfigure}[h]{0.44\textwidth}
        \centering
        \includegraphics[width=\textwidth]{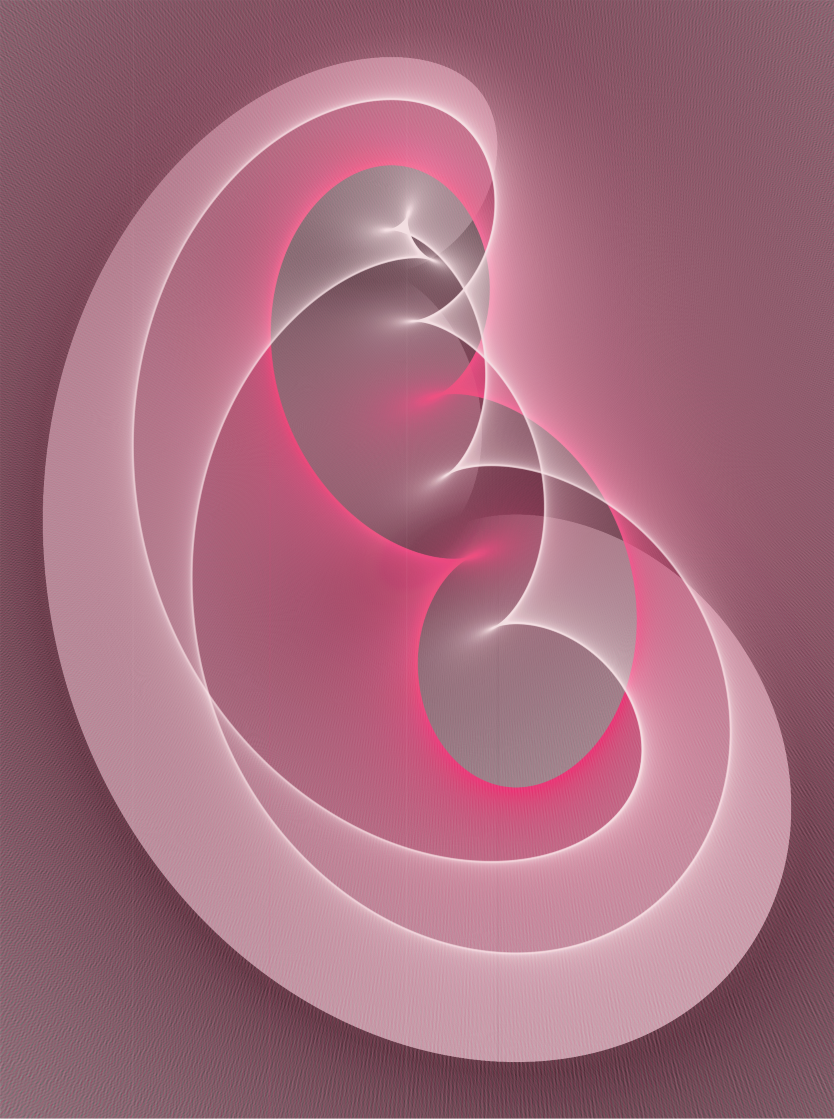}
        \caption{\textit{The Beginning}}
        \label{fig: begining}
    \end{subfigure}
    \caption{Artistic visualizations in \texttt{Mathematica} (part 1)}
    \label{fig:finalexamples1}
\end{figure}

\begin{figure}[h]
    \centering
    \begin{subfigure}[h]{0.44\textwidth}
        \centering
        \includegraphics[width=\textwidth]{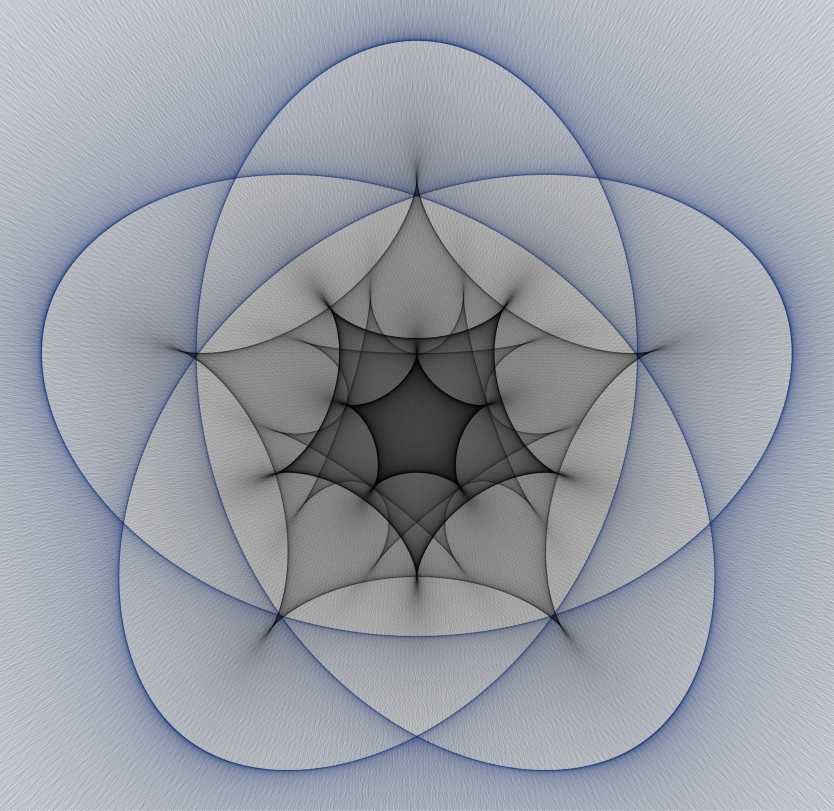}
        \caption{\textit{MiNI's Logo}}
        \label{fig: logomini}
    \end{subfigure}
    \hfill
    \begin{subfigure}[h]{0.44\textwidth}
        \centering
        \includegraphics[width=\textwidth]{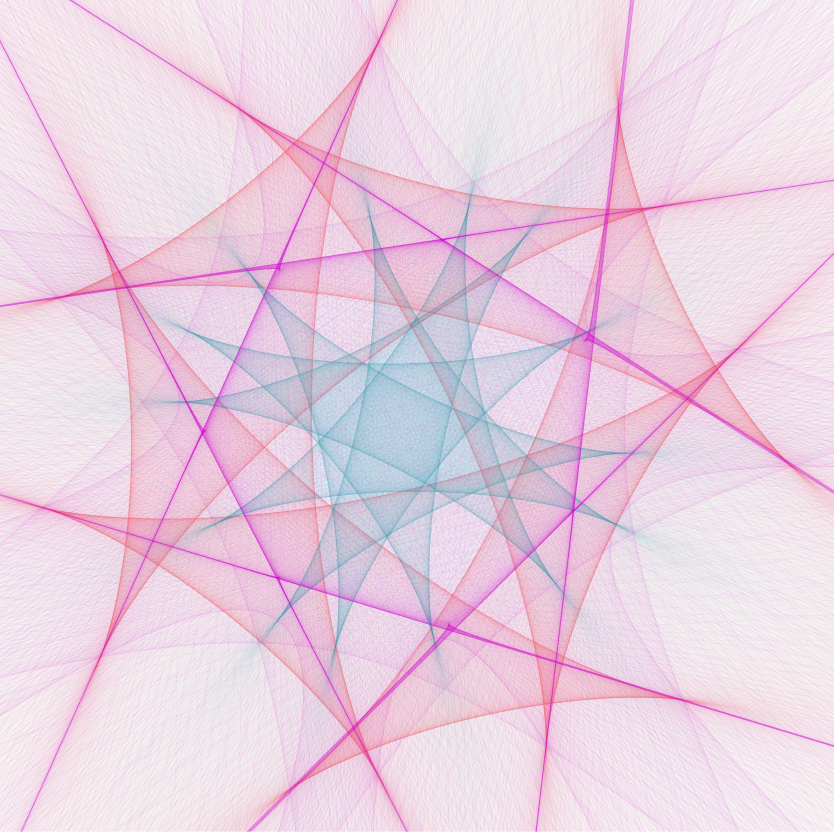}
        \caption{\textit{Chaos}}
        \label{fig: asymmetry}
    \end{subfigure}
    \\ 
    \begin{subfigure}[h]{0.44\textwidth}
        \centering
        \includegraphics[width=\textwidth]{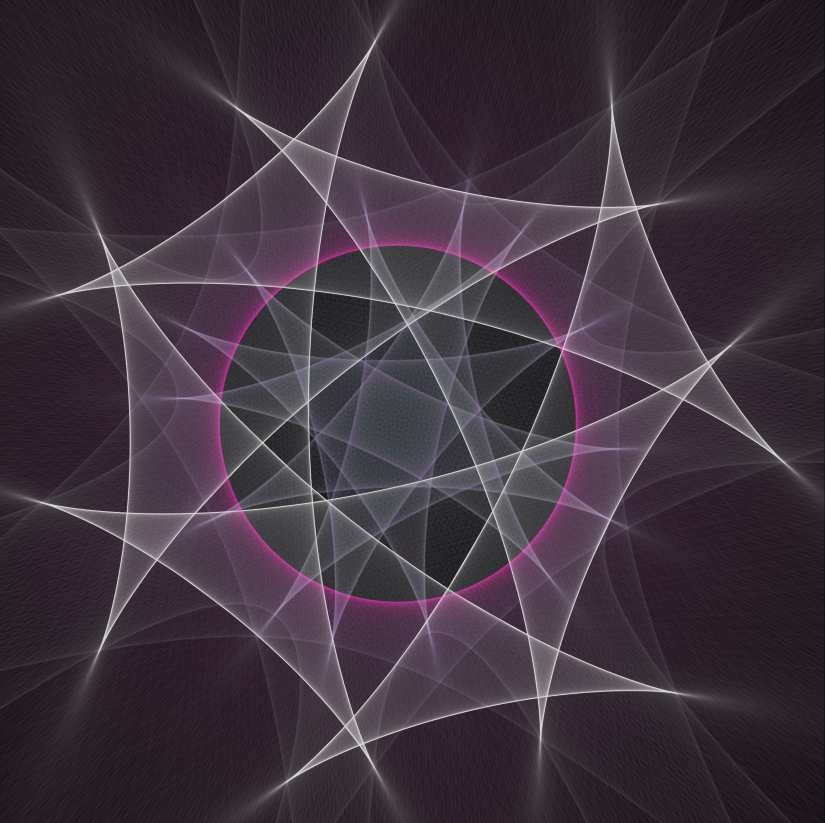}
        \caption{\textit{Shuriken}}
        \label{fig: shuriken}
    \end{subfigure}
    \hfill
    \begin{subfigure}[h]{0.44\textwidth}
        \centering
        \includegraphics[width=\textwidth]{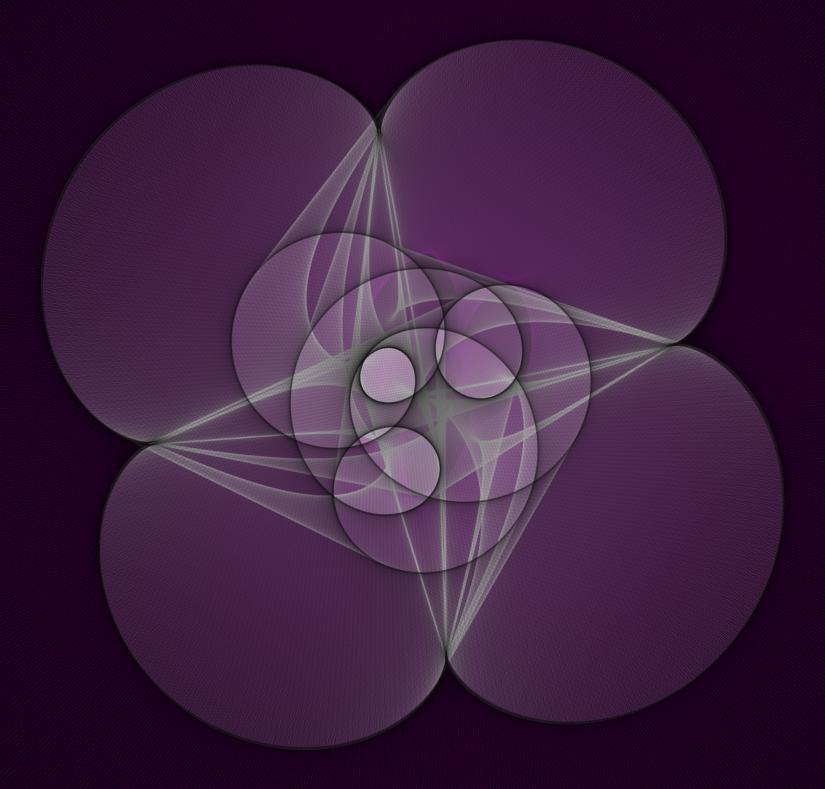}
        \caption{\textit{Singular Four-Leaf Clover}}
        \label{fig: twocurves}
    \end{subfigure}
    \label{fig:dpsz}
    \begin{subfigure}[h]{0.44\textwidth}
        \centering
        \includegraphics[width=\textwidth]{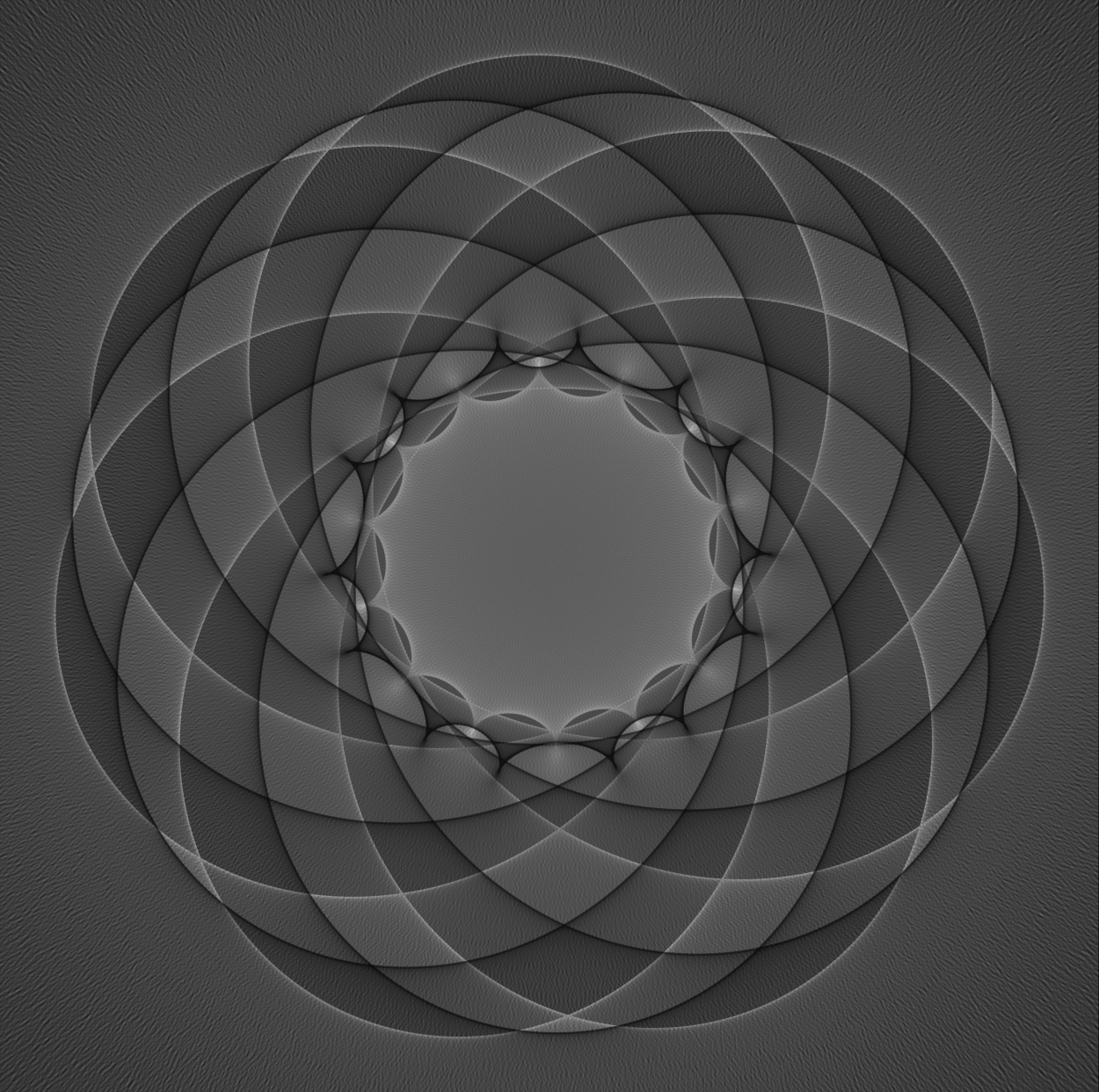}
        \caption{\textit{Black and White}}
        \label{fig: blackandwhite}
    \end{subfigure}
    \hfill
    \begin{subfigure}[h]{0.44\textwidth}
        \centering
        \includegraphics[width=\textwidth]{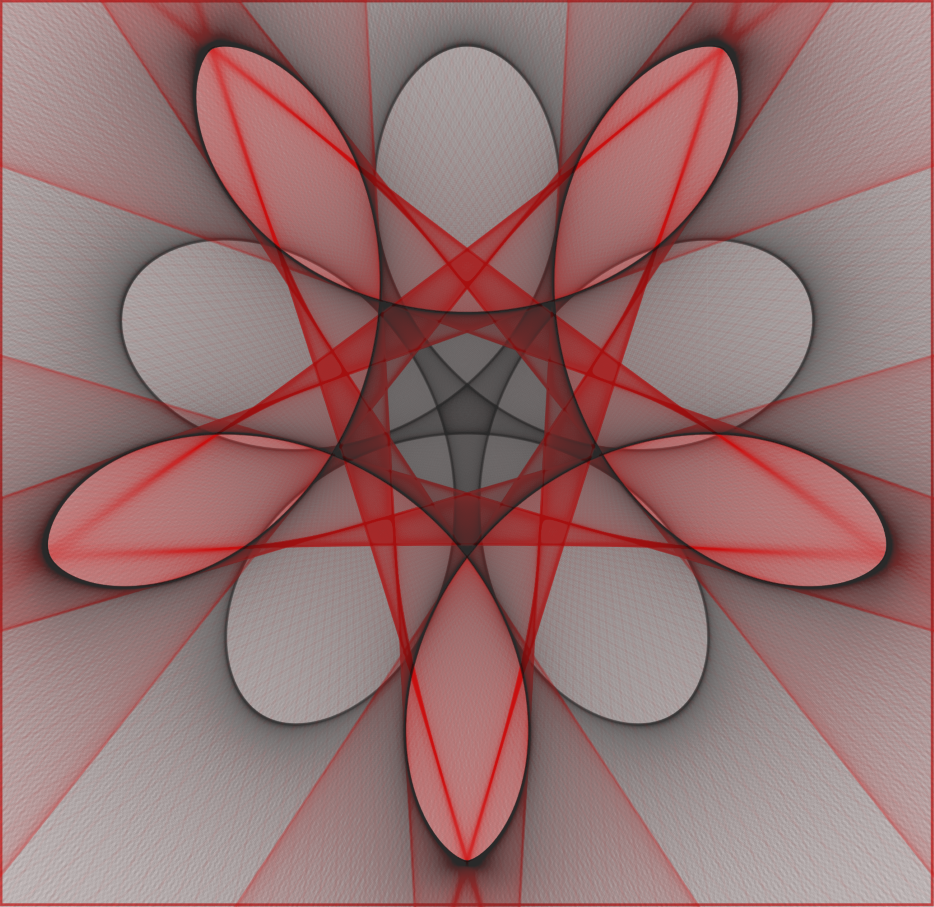}
        \caption{\textit{Beast's Rosette}}
        \label{fig: beastrosette}
    \end{subfigure}
    \caption{Artistic visualizations in \texttt{Mathematica} (part 2)}
    \label{fig:finalexamples2}
\end{figure}

\begin{figure}[h]
    \centering
    \begin{subfigure}[h]{0.44\textwidth}
        \centering
        \includegraphics[width=\textwidth]{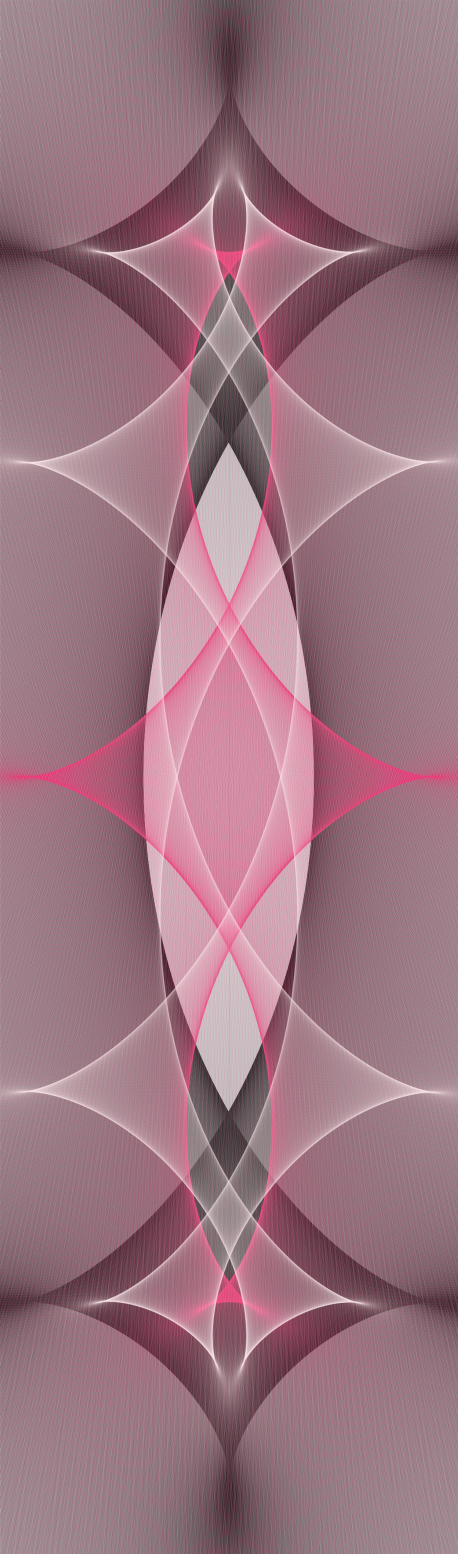}
        \caption{\textit{Harmony}}
        \label{fig: harmony}
    \end{subfigure}
    $\ \ $
    \begin{subfigure}[h]{0.475\textwidth}
        \centering
        \includegraphics[width=\textwidth]{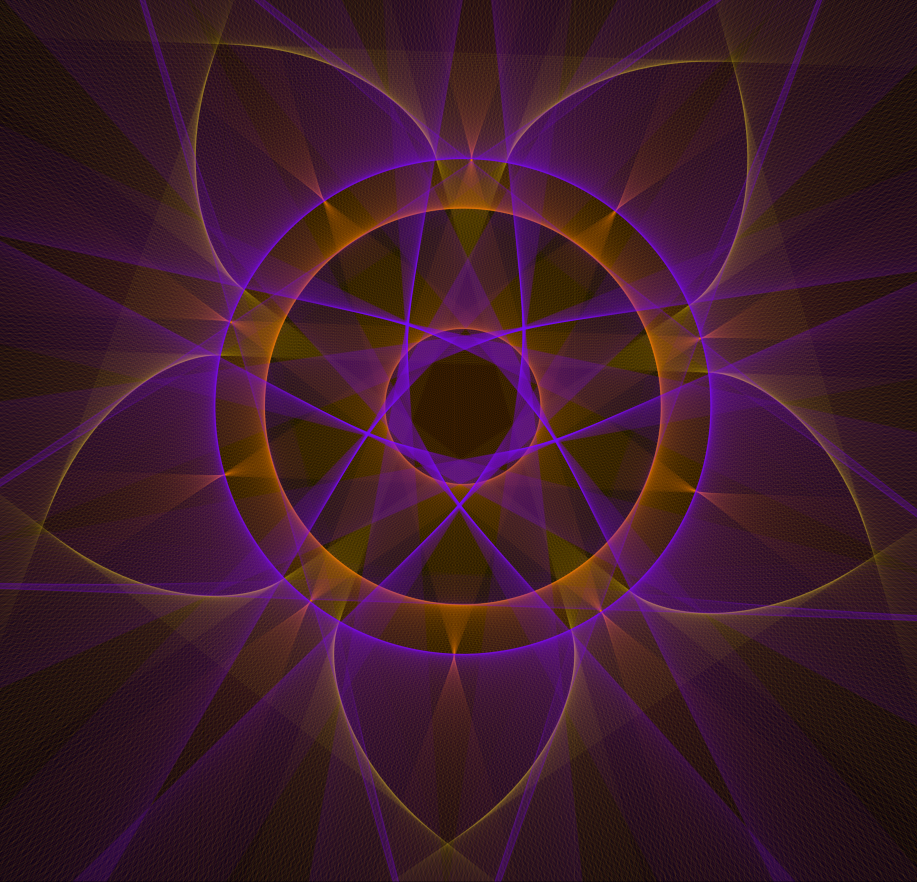}
        \caption{\textit{L\'{o}th}}
        \label{fig: harmony}
        \includegraphics[width=\textwidth]{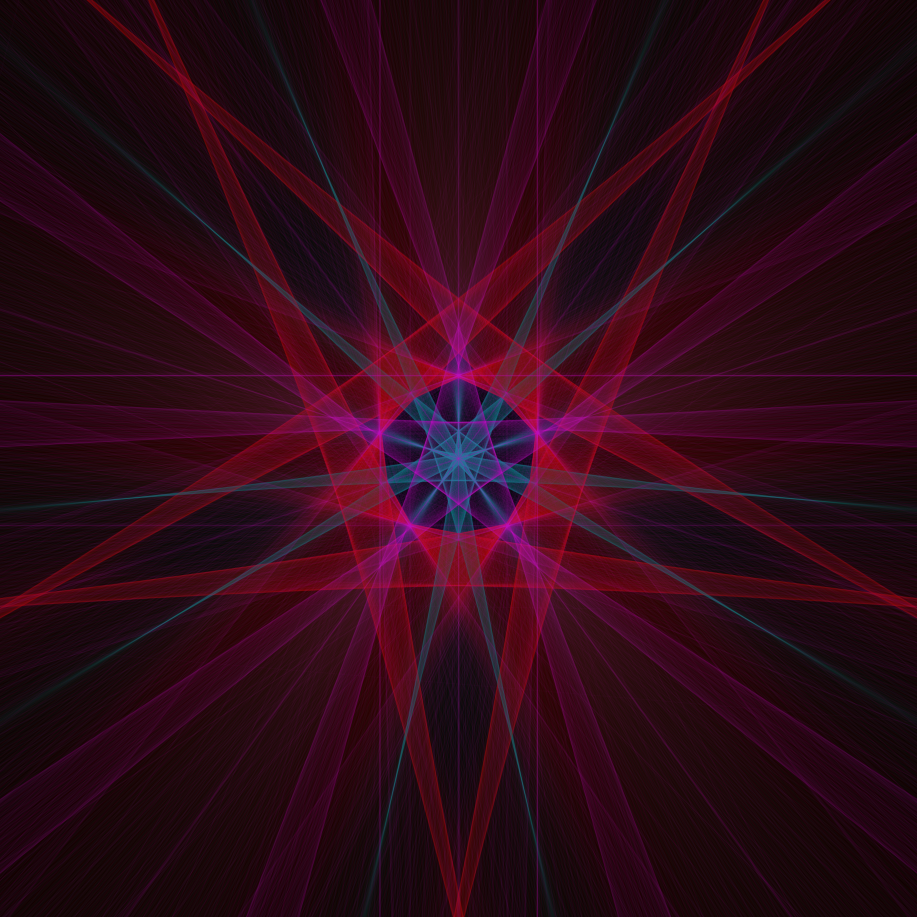}
        \caption{\textit{Pentaception}}
        \label{fig: harmony}
        \includegraphics[width=\textwidth]{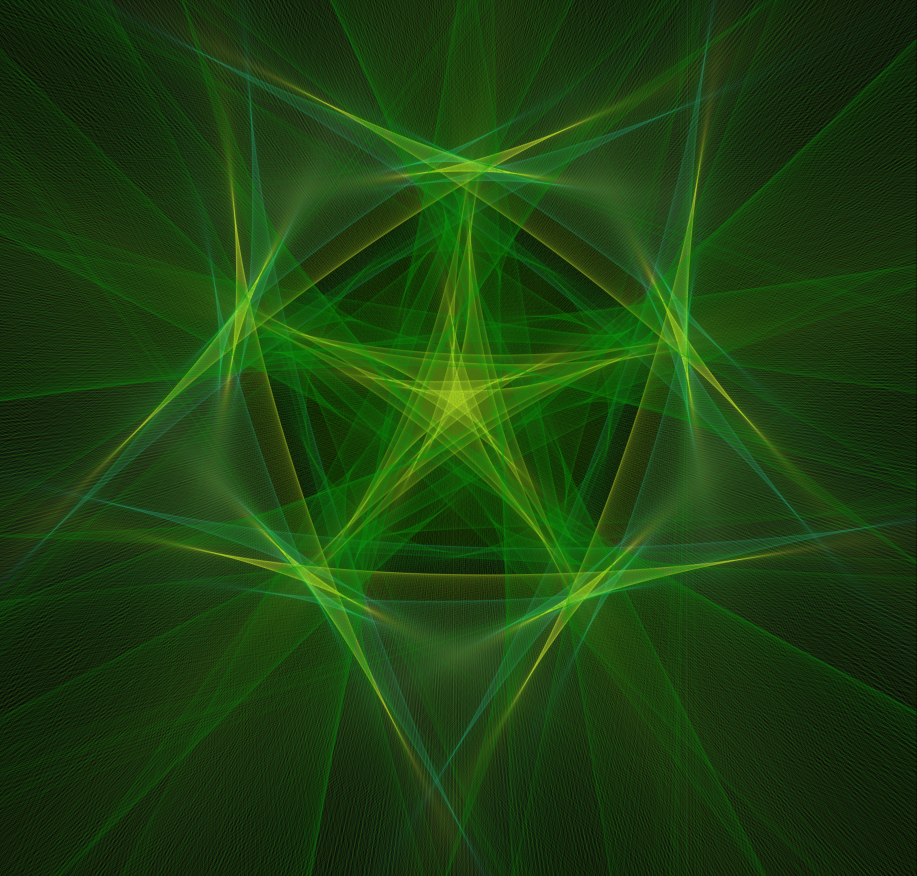}
        \caption{\textit{Pestilence}}
        \label{fig: harmony}
    \end{subfigure}
    \caption{Artistic visualizations in \texttt{Mathematica} (part 3)}
    \label{fig:finalexamples3}
\end{figure}

This section unveils the drawings we collaboratively created, based on the techniques discussed in the previous chapter. We are thrilled to showcase these artworks, and we fervently hope they will resonate with a diverse audience. We invite you to explore and enjoy them (Figures \ref{fig:finalexamples1}, \ref{fig:finalexamples2}, and \ref{fig:finalexamples3}). 

\clearpage
\pagebreak

\clearpage
\pagebreak

\clearpage
\pagebreak

\clearpage
\pagebreak

\section*{Acknowledgements}
The Authors would like to thank Professor Stanis\l{}aw Janeczko for the helpful discussions, Center for Advanced Studies at Warsaw University of Technology, 
and Mathematics Popularization Day for the opportunity to exhibit their works to a~wider audience.

\bibliographystyle{amsalpha}

\end{document}